\documentclass[10pt,a4paper]{article}
\usepackage{amsfonts}
\usepackage{amssymb}
\usepackage{mathrsfs}
\usepackage{amsthm}
\usepackage{setspace}
\usepackage{color}

\usepackage{fancyhdr}
\newtheorem{theo}{Theorem}[section]
\newtheorem{prop}[theo]{Proposition}
\newtheorem{cor}[theo]{Corollary}
\newtheorem{lemm}[theo]{Lemma}

\theoremstyle{definition}
\newtheorem*{rem}{Remark}
\newtheorem{defi}[theo]{Definition}

\newenvironment{lproof}{\emph{Proof of Lemma.}}{ \qed \par}

\newcommand{\be}{\begin{eqnarray*}}
\newcommand{\ee}{\end{eqnarray*}}
\newcommand{\beqa}{\begin{eqnarray}}
\newcommand{\eeqa}{\end{eqnarray}}
\newcommand{\ba}{\begin{array}}
\newcommand{\ea}{\end{array}}
\newcommand{\onab}{\overrightarrow{\nabla}}
\newcommand{\mc}{\mathcal}
\newcommand{\mf}{\mathfrak}

\newcommand{\rP}{\mathsf{P}}

\newcommand{\mbb}{\mathbb}
\newcommand{\wt}{\widetilde}
\newcommand{\wh}{\widehat}
\newcommand{\A}{\mc{A}}
\newcommand{\B}{\mc{B}}
\newcommand{\T}{\mc{T}}

\newcommand\CR{\begin{picture}(46,12)%
\put(5,3){\line(1,0){16}}%
\put(25,3){\line(1,0){16}}%
\put(3,3){\makebox(0,0){$\times$}}%
\put(23,3){\makebox(0,0){$\circ$}}%
\put(43,3){\makebox(0,0){$\times$}}%
\end{picture}}

\newcommand\Bpalgebra{\begin{picture}(76,12)\put(5,3){\line(1,0){16}}%
\put(25,3){\line(1,0){6}}%
\put(51,3){\line(-1,0){6}}\put(55,4){\line(1,0){17}}%
\put(54,2){\line(1,0){18}}\put(39,3){\makebox(0,0){\dots}}%
\put(63,2){\makebox(0,0){$>$}}%
\put(3,3){\makebox(0,0){$\circ$}}%
\put(23,3){\makebox(0,0){$\circ$}}%
\put(53,3){\makebox(0,0){$\circ$}}%
\put(73,3){\makebox(0,0){$\times$}}%
\end{picture}}

\newcommand\Dpalgebra{\begin{picture}(76,12)\put(5,3){\line(1,0){16}}%
\put(25,3){\line(1,0){6}}%
\put(51,3){\line(-1,0){6}}%
\put(55,4){\line(2,1){17}}%
\put(54,2){\line(2,-1){18}}%
\put(39,3){\makebox(0,0){\dots}}%
\put(3,3){\makebox(0,0){$\circ$}}%
\put(23,3){\makebox(0,0){$\circ$}}%
\put(53,3){\makebox(0,0){$\circ$}}%
\put(74,12){\makebox(0,0){$\circ$}}%
\put(73,-7){\makebox(0,0){$\times$}}%
\end{picture}}

\newcommand\Dpalgebrajaw{\begin{picture}(76,12)\put(5,3){\line(1,0){16}}%
\put(25,3){\line(1,0){6}}%
\put(51,3){\line(-1,0){6}}%
\put(55,4){\line(2,1){17}}%
\put(54,2){\line(2,-1){18}}%
\put(39,3){\makebox(0,0){\dots}}%
\put(3,3){\makebox(0,0){$\circ$}}%
\put(23,3){\makebox(0,0){$\circ$}}%
\put(53,3){\makebox(0,0){$\circ$}}%
\put(74,12){\makebox(0,0){$\times$}}%
\put(73,-7){\makebox(0,0){$\times$}}%
\end{picture}}

\newcommand\Balgebra{\begin{picture}(76,12)\put(5,3){\line(1,0){16}}%
\put(25,3){\line(1,0){6}}%
\put(51,3){\line(-1,0){6}}\put(55,4){\line(1,0){17}}%
\put(54,2){\line(1,0){18}}\put(39,3){\makebox(0,0){\dots}}%
\put(63,2){\makebox(0,0){$>$}}%
\put(3,3){\makebox(0,0){$\circ$}}\put(23,3){\makebox(0,0){$\circ$}}%
\put(53,3){\makebox(0,0){$\circ$}}\put(73,3){\makebox(0,0){$\circ$}}%
\end{picture}}

\newcommand\Dalgebra{\begin{picture}(76,12)\put(5,3){\line(1,0){16}}%
\put(25,3){\line(1,0){6}}%
\put(51,3){\line(-1,0){6}}%
\put(55,4){\line(2,1){17}}%
\put(54,2){\line(2,-1){18}}%
\put(39,3){\makebox(0,0){\dots}}%
\put(3,3){\makebox(0,0){$\circ$}}\put(23,3){\makebox(0,0){$\circ$}}%
\put(53,3){\makebox(0,0){$\circ$}}%
\put(74,12){\makebox(0,0){$\circ$}}%
\put(74,-7){\makebox(0,0){$\circ$}}%
\end{picture}}

\newcommand\trione{\begin{picture}(49,12)\put(24,3){\line(-1,0){18}}%
\put(28,4){\line(2,1){17}}%
\put(27,2){\line(2,-1){18}}%
\put(4,3){\makebox(0,0){$\circ$}}%
\put(26,3){\makebox(0,0){$\circ$}}%
\put(47,12){\makebox(0,0){$\circ$}}%
\put(46,-7){\makebox(0,0){$\times$}}%
\end{picture}}

\newcommand\tritwo{\begin{picture}(49,12)\put(24,3){\line(-1,0){18}}%
\put(28,4){\line(2,1){17}}%
\put(27,2){\line(2,-1){18}}%
\put(5,2){\makebox(0,0){$\times$}}%
\put(26,3){\makebox(0,0){$\circ$}}%
\put(47,12){\makebox(0,0){$\circ$}}%
\put(47,-7){\makebox(0,0){$\circ$}}%
\end{picture}}

\begin{document}

\title{Free $3$-distributions: holonomy, Fefferman constructions and dual distributions}
\author{Stuart Armstrong, \\ Fakult\"at f\"ur Mathematik, Universit\"at Wien, Nordbergstr. 15, 1090 Wien, Austria}
\date{2008}
\maketitle

\begin{abstract}
This paper analyses the parabolic geometries generated by a free $3$-distribution in the tangent space of a manifold. It shows the existence of normal Fefferman constructions over CR and Lagrangian contact structures corresponding to holonomy reductions to $SO(4,2)$ and $SO(3,3)$, respectively. There is also a fascinating construction of a `dual' distribution when the holonomy reduces to $G_2'$. The paper concludes with some holonomy constructions for free $n$-distributions for $n>3$.
\end{abstract}

\section{Introduction}

On a manifold $M$, let $H \subset TM$ be a distribution of rank $n$. Then there is a well defined map $\mc{L}: H \wedge H \to TM/H$. For $X,Y$ sections of $H$, it is given by the quotiented Lie bracket $X\wedge Y \to [X,Y]/H$. Then $H$ is a \emph{free} $n$-distribution if $\mc{L}$ is an isomorphism. The moniker ``free'' comes from the fact that there are no relations between sections of $H$ that would cause $\mc{L}$ to fail injectivity.

This condition immediately implies that $TM/H$ is of rank $n(n-1)/2$, thus that $M$ is of dimension $m = n(n+1)/2$. Bryant \cite{bryskew} has achieved some major results in the case of $n=3, m=3$, a free $3$-distribution in a $6$-manifold, but a lot remains unknown.

Fortunately, these structures lead themselves to be treated with the general tools of Cartan connections on parabolic geometries (\cite{TCPG} and \cite{CartEquiv}). The homogeneous model is provided by the set of maximal isotropic planes in $\mathbb{R}^{n+1,n}$. The group of transformations is $G = PSO(n+1,n)$ (with Lie algebra $\mf{g}$), while the stabiliser of a point is $P = GL(n) \rtimes \mathbb{R}^n \rtimes \wedge^2 \mathbb{R}^n$. The stabiliser´s Lie algebra is $\mf{p}$ which has nilradical $\mathbb{R}^n \rtimes \wedge^2 \mathbb{R}^n$. These are precisely the two-step free nilpotent Lie algebras, with the Lie bracket from $\mathbb{R}^n \otimes \mathbb{R}^n$ to $\wedge^2 \mathbb{R}^n$ being given by taking the wedge. Using the Killing form on $\mf{g}$, we may define $\mf{g}_{(-1)}$ as $[\mf{g},\mf{p}^{\perp}]$; by the action of $G$, this vector space gives a subbundle of $TG$, that then projects down to a distribution on $G/P$ of rank $n$. This distribution is precisely the free $n$-distribution for the homogeneous model.

In the general (non-homogeneous) case, we do not introduce any extra information, or make any choices upon taking the Cartan connection, as the normal Cartan connection for a free $n$-distribution is determined entirely by $H$ (\cite{two}).

The $n=3, m=6$ case is the most tractable (as the normal Cartan connection is torsion-free), and in many ways the most interesting. The free $3$-distribution has a Fefferman construction into a conformal structure on the same manifold \cite{bryskew} (a special case of the Fefferman construction form free $n$-distributions into almost spinorial structures). We will show that the Fefferman construction is is normal if the connection coming from the free $3$-distribution is normal. Conversely, if the holonomy of a normal conformal Cartan connection reduces to $Spin(3,4)$, it must locally the Fefferman construction of a normal Cartan connection of a free $3$-distribution, meaning that we have many known examples of holonomy reductions \cite{mecon}.

Using this, we can then show that holonomy reductions exist to the groups $G_2'$, $SU(2,2) \cong Spin(4,2)_0$ and $SL(4,\mbb{R}) \cong Spin(3,3)_0$. The reduction to $G_2'$ does not come from any Fefferman construction, but has a fascinating geometry. On an open dense set of the manifold, this holonomy reduction generates a canonical Weyl structure $\nabla$. This $\nabla$ determines a splitting of $T = T_{-2} \oplus H$, where $H$ is the free $3$-distribution defining the Cartan connection. Then $H' = T_{-2}$ is also a free $3$-distribution, and the normal Tractor connection it generates is isomorphic with $\onab$. Under certain bundle inclusions, the two Tractor connections are the same. Iterating this procedure with $H'$ generates $H$ again; thus $H$ and $H'$ are in some sense `dual' distributions.

The other two holonomy reductions, $SU(2,2)$ and $SL(4,\mbb{R})$ arise from their own Fefferman constructions -- over integrable CR manifolds and integrable Lagrangian contact structures, respectively. Here, normality of the underlying Cartan connections is equivalent with normality of the one generated by the free $3$-distribution.

Finally, the paper concludes in an appendix, with a look at some simple holonomy algebras in the higher rank cases where $n > 3$. It demonstrates that there exists $n$-distributions with holonomy algebras isomorphic with $\mbb{R}^p$ for all $p<n(n-1)/2$.

\subsection*{Acknowledgements}
It gives me great pleasure to acknowledge the financial support of an ESI Junior Fellowship program and project P19500-N13 of the ``Fonds zur F\"orderung der wissenschaftlichen Forschung (FWF)'', as well as the help, proofreading and comments of Andreas {\v{C}}ap and Daniel Fox.

\section{The geometry of free $n$-distributions}

\subsection{Homogeneous model} \label{hom:mod}

The homogeneous model for a free $n$-distribution is the space of isotropic $n$-planes inside $\mathbb{R}^{(n+1,n)}$. Let $\mf{g} = \mf{so}(n+1,n)$, and put the metric on $\mathbb{R}^{(n+1,n)}$ in the form:
\be
\left( \begin{array}{ccc} 0&0&Id_n \\ 0&1&0 \\ Id_n &0&0 \end{array} \right),
\ee
where $Id_n$ is the identity matrix on $\mathbb{R}^n$. The algebra $\mf{g}$ is then spanned by elements of the form:
\be
\left( \begin{array}{ccc} A&v&B \\ w&0&-v^t \\ C&-w^t&-A^t \end{array} \right),
\ee
where $A \in \mf{gl}(n)$, $B^t = - B$ and $C^t = -C$. Then the isotropic $n$-plane
\be
V = \{(a_1, a_2, \ldots, a_n, 0, 0, \ldots, 0) | a_j \in \mathbb{R}\}
\ee
is preserved by the subalgebra $\mf{p} \subset \mf{g}$ consisting of elements with $C = 0$ and $w = 0$. This algebra is isomorphic to $\mf{gl}(n) \oplus (\mathbb{R}^n)^* \oplus (\wedge^2 \mathbb{R}^n)^*$, with the natural algebraic structure.

The homogeneous manifold is $M = G/P$, where $G = PSO(n+1,n)_0$ and $P \subset G$ is the Lie group with Lie algebra $\mf{p}$. $M$ is of dimension $(2n+1)(2n)/2 - (n^2 + n + n(n-1)/2) = n(n+1)/2$. There is a subspace $\mf{g}_{(-1)}$ of $\mf{g}$, consisting of those elements with $C = 0$. Left action on $G$ generalises this to a distribution $\widehat{H} \subset TG$. This distribution is preserved by $P$, and hence the map $G \to M$ maps $\widehat{H}$ to a distribution $H \subset TM$. This distribution is evidently of rank $n$.

Now consider the differential bracket of left-invariant vector fields which are sections of $\widehat{H}$; this matches up with the Lie bracket on $\mf{g}_{(-1)}$, thus $[\widehat{H}, \widehat{H}]$ spans all of $TG$ (by $[\widehat{H}, \widehat{H}]$ we mean the bundle spanned by the bracket of all pairs of sections of $\wh{H}$). Consequently, $[H,H]$ must span all of $TM$, making $H$ free by dimensional considerations.

\subsection{Cartan connection}

Given a semi-simple Lie algebra $\mathfrak{g}$ with Killing form $(-,-)$, a subalgebra $\mf{p} \subset \mf{g}$ is said to be \emph{parabolic} iff $\mf{p}^{\perp}$ is the \emph{nilradical} of $\mf{p}$, i.e.~its maximal nilpotent ideal. This gives (\cite{paradef}, details are also available in the author's thesis \cite{methesis}) a filtration of $\mf{g}$:
\be
\begin{array}{ccccccccccccc}
\mf{g}_{(-k)}& \supset & \mf{g}_{(k-1)} & \supset &\ldots& \supset &\mf{g}_{(0)}& \supset &\mf{g}_{(1)}& \supset & \ldots & \supset & \mf{g}_{(k)}, \\
| | &&&&&&| | && | | && \\
\mf{g} &&&&&&\mf{p}&& \mf{p}^{\perp}&&
\end{array}
\ee
such that $\{ \mf{g}_{(j)} , \mf{g}_{(l)} \} \subset \mf{g}_{(j+l)}$. The associated graded algebra is $ {\textrm{gr}}(\mf{g}) = \bigoplus_{-k}^{k} \mf{g}_j$, where $\mf{g}_j = \mf{g}_{(j)}/\mf{g}_{(j+1)}$. By results from \cite{paradef}, ${\textrm{gr}}(\mf{g})$ is isomorphic to $\mf{g}$. Furthermore, there is a unique element $\epsilon_0$ in $\mf{g}_0$ such that $\{ \epsilon_0, \xi \} = j \xi$ for all $\xi \in \mf{g}_j$. The isomorphisms ${\textrm{gr}}(\mf{g}) \cong \mf{g}$ compatible with the filtration are then given by a choice of lift $\epsilon$ of $\epsilon_0$ with respect to the exact sequence
\beqa \label{p:sequence}
0\to \mf{p}^{\perp} \to \mf{p} \to \mf{g}_0 \to 0.
\eeqa
This means that the gradings of $\mf{g}$ compatible with the filtration form an affine space modelled on the nilradical $\mf{p}^{\perp}$. Define $P$ and $G_0$ as the subgroups of $G$ that preserve the filtration and the grading, respectively. It is easy to see that their Lie algebras are $\mf{p}$ and $\mf{g}_0$, and that the inclusion $G_0 \subset G$ is non-canonical.

\begin{defi}[Cartan connection]
A Cartan connection on a manifold $M$ for the parabolic subalgebra $\mf{p} \subset \mf{g}$ is given by a principal $P$-bundle
\be
\mc{P} \to M,
\ee
and a one form $\omega \in \Omega^1(\mc{P}, \mf{g})$ with values in the Lie algebra $\mf{g}$ such that:
\begin{enumerate}
\item $\omega$ is equivariant under the $P$-action ($P$ acting by $Ad$ on $\mathfrak{g}$),
\item $\omega(\sigma_A) = A$, where $\sigma_A$ is the fundamental vector field of $A \in \mathfrak{p}$,
\item $\omega_u: T \mathcal{P}_u \to \mathfrak{g}$ is a linear isomorphism for all $u \in \mathcal{P}$.
\end{enumerate}
\end{defi}

The inclusion $P \subset G$ generates a bundle inclusion $\mc{P} \subset \mc{G}$ and $\omega$ is then the pull back of a unique $G$-equivariant connection form on $\mc{G}$ which we will also designate by $\omega$. Since $\omega$ takes values in $\mf{g}$, it generates a standard connection on any vector bundle associated to $\mc{G}$. This connection is called the \emph{Tractor connection} and will be designated by $\onab$.

In the case we are looking at, $\mf{g} = \mf{so}(n+1,n)$ and $\mf{p} = \mf{gl}(n) \oplus (\mathbb{R}^n)^* \oplus (\wedge^2 \mathbb{R}^n)^*$. $\mf{g}_0$ is simply $\mf{gl}(n)$, and the nilradical of $\mf{p}$ is $\mf{p}^{\perp} = (\mathbb{R}^n)^* \oplus (\wedge^2 \mathbb{R}^n)^*$. In terms of the notation for parabolic subalgebras introduced in \cite{Beastwood}, this is:
\Bpalgebra. Now define the Lie algebra bundle
\be
\mc{A} = \mc{P} \times_P \mf{g}.
\ee
This inherits an algebra\"ic bracket $\{,\}$ coming from the Lie bracket $[,]$ on $\mf{g}$, and has a natural filtration
\be
\mc{A} =\mc{A}_{(-2)} \supset \mc{A}_{(-1)} \supset \mc{A}_{(0)} \supset \mc{A}_{(1)} \supset \mc{A}_{(2)},
\ee
with $\mc{A}_{(j)} = \mc{P} \times_P \mf{g}_{(j)}$. Paper \cite{TCPG} demonstrates that the tangent space $T$ of $M$ is equal to the quotient bundle $\mc{A} / \mc{A}_{(0)}$. The killing form gives an isomorphism
\be
\mc{A}_{(-1)} \cong (\mc{A} / \mc{A}_{(0)})^* = T^*.
\ee
Hence there is a well defined inclusion $T^* \subset \mc{A}$, and a well defined projection $\mc{A} \to T$. We define the \emph{graded} bundles $\A_j$ as:
\be
\mc{A}_j = \mc{P} \times_P \mf{g}_j = \A_{(j)} / \A_{(j+1)}.
\ee

\begin{defi}[Weyl structure]
A \emph{Weyl structure} on $(M, \mc{P}, \omega)$ is a $P$-equivariant function $\eta: \mc{P} \to \mathfrak{p}$ that is always a lift of the grading element $\epsilon_0$, as in equation (\ref{p:sequence}).
\end{defi}
A Weyl structure gives a splitting of $\mf{g}$, and consequently allows a decomposition of both the Lie algebra bundle and the Cartan connection:
\be
\mc{A} &=& \mc{A}_{-2} \oplus \mc{A}_{-1} \oplus \mc{A}_{0} \oplus \mc{A}_{1} \oplus \mc{A}_{2} \\
\omega &=& \omega_{-2} + \omega_{-1} + \omega_0 + \omega_1 + \omega_2.
\ee
Dividing out by the action of $\mf{p}^{\perp}$ gives a quotient map $\mf{p} \to \mf{g}_0$ and hence a bundle map $\mc{P} \to \mc{G}_0$. There is a unique $G_0$-equivariant one-form on $\mc{G}_0$ of which $\omega_0$ is the pull-back; we will call it $\omega_0$ as well. This is a principal connection on $\mc{G}_0$; since $G_0$ acts on $\mf{g}_{-2} + \mf{g}_{-1} \cong \mf{g}/\mf{p}$, then $T = \mc{A} / \mc{A}_{(0)}$ is an associated bundle to $\mc{G}_0$. This implies that $\omega_0$ generates an affine connection $\nabla$ on the tangent bundle.

These $\nabla$'s are called preferred connections; they are in one-to-one correspondence with Weyl structures and hence to compatible splittings of $\mc{A}$. They consequently form an affine space modelled on $\mc{A}_{(-1)} = T^*$; the relation between two preferred connections $\nabla$ and $\widehat{\nabla}$ is given explicitly by the one-form $\Upsilon$ on $M$ such that
\beqa \label{change:con}
\nabla_X Y = \widehat{\nabla}_X Y + \{\{X, \Upsilon \}, Y\}_{-},
\eeqa
with $\{,\}$ the Lie bracket on $\mc{A}$ and using the natural inclusion $T = \A/\A_{(0)} \cong \A_{-2} \oplus \A_{-1} \subset \A$ given by the splitting. The splitting of $\mc{A}$ gives further splittings $T = T_{-2} \oplus T_{-1}$ and $T^* = T^*_{1} \oplus T^*_{2}$. The bundles $T_{-1}$ and $T^*_{2} = (T_{-1})^{\perp}$ are defined independently of the Weyl structure, since they are preserved by the action of $T^* \subset \mc{A}$.

Given a preferred $\nabla$, the Tractor connection on $\mc{A}$ and any associated bundles is given by:
\be
\onab_X v = \nabla_X v + X \cdot v + \rP(X) \cdot v,
\ee
where $\rP$ is the rho-tensor, a section of $T^* \otimes T^*$, and $\cdot$ is the action of $T$ and $T^*$ on $\mc{A}$ given by the Lie bracket. There are various sign conventions for $\rP$ (especially in conformal geometry); in this paper, we will take the sign convention that makes the above formula true.

\begin{defi}[Curvature]
The curvature of the Cartan connection is defined to be the two-form $\kappa = d \omega + \frac{1}{2}\{ \omega , \omega \} \in \Omega^2(\mc{P}, \mf{g})$. It is easy to see that $\kappa$ vanishes on vertical vectors, and is $P$-equivariant; consequently dividing out by the action of $P$, $\kappa$ may seen as an element of $\Omega^2(M, \mathcal{A})$; in this setting, it is the curvature of the Tractor connection $\onab$. Finally, the inclusion $T^* \subset \mc{A}$ implies that $\kappa$ is equivalent to a $P$-equivariant function from $\mc{P}$ to $\wedge^2 \mf{p}^{\perp} \otimes \mf{g}$. We shall use the designation $\kappa$ interchangeably for these three equivalent definitions, though it is generally the third one we shall be using.
\end{defi}

Given a grading on $\mf{g}$, there is a decomposition of any tensor product $\otimes^c \mf{g} = \sum \mf{g}_{j_1, j_2, \ldots j_c} $ where $\mf{g}_{j_1, j_2, \ldots j_c} = \mf{g}_{j_1} \otimes \mf{g}_{j_2} \otimes \ldots \otimes \mf{g}_{j_c}$. The homogeneity of $\mf{g}_{j_1, j_2, \ldots j_c}$ is defined to be the sum $j_1 + j_2 + \ldots + j_c$. Any element of $\eta$ of $\otimes^c \mf{g}$ can be decomposed into homogeneous elements $\eta_{j_1, j_2, \ldots, j_c}$. The \emph{minimal homogeneity} of $\eta$ is defined to be the lowest homogeneity among the non-zero $\eta_{j_1, j_2, \ldots, j_c}$.

Homogeneity is not preserved by the action of $P$; however since $\mf{p}$ consists of elements of homogeneity zero and above, the minimal homogeneity of any element is preserved by the action of $P$. Since $\kappa$ is a map to $\wedge^2 \mf{p}^{\perp} \otimes \mf{g}$, the following definition makes sense:
\begin{defi}[Regularity]
A Cartan connection is regular iff the minimal homogeneity of $\kappa$ is greater than zero.
\end{defi}

There are well defined Lie algebra differentials $\partial: \wedge^c \mathfrak{p}^{\perp} \otimes \mf{g} \to \wedge^{c+1} \mathfrak{p}^{\perp} \otimes \mf{g}$ and codifferentials $\partial^*: \wedge^c \mathfrak{p}^{\perp} \otimes \mf{g} \to \wedge^{c-1} \mathfrak{p}^{\perp} \otimes \mf{g}$. In terms of decomposable elements, the co-differential is given by
\be
\partial^* (u_1 \wedge \ldots \wedge u_c \otimes v) &=& \sum_{j <k} (-1)^{j+k} \{u_j, u_k \} \wedge u_1 \wedge \ldots \wedge u_c \otimes v + \\
&& \sum_{j} (-1)^{j+1} u_1 \wedge \ldots  \wh{u_j} \ldots \wedge u_c \otimes \{u_j, v\}.
\ee
\begin{defi}[Normality] \label{normal:def}
A Cartan connection is normal iff $\partial^* \kappa = 0$.
\end{defi}
If we let $(Z_l)$ be a frame for $T$ and $(Z^l)$ a dual frame for $T^*$, this condition is given, in terms of $\kappa$ an element of $\Omega^2(M, \mc{A})$, as
\beqa \label{normal:form}
(\partial^* \kappa)(X) = \sum_{l} \{ Z^l, \kappa_{(Z_l, X)} \} - \frac{1}{2} \kappa_{(\{Z^l,X \}_{-}, Z_l)},
\eeqa
for all $X \in \Gamma(T)$.

For free $n$-distribution, we have $G_0 = GL(\mf{g}_{-1})$. Paper \cite{CartEquiv} then demonstrates that a regular, normal Cartan connection for these Lie groups is determined entirely by the distribution $T_{-1}$. We shall call this distribution $H$, as before. Since $\omega$ is regular, the algebraic bracket $\{ ,\}$ matches up with graded version of the Lie bracket $[,]$ of vector fields; in other words, if $X$ and $Y$ are sections of $H$,
\be
[X,Y]/H - \{X,Y\}/H = 0,
\ee
implying (since $\{,\}$ is surjective $H\wedge H \to TM/H$) that $H$ is maximally non-integrable. Indeed, any maximally non-integrable $H$ of correct dimension and co-dimension determines a unique normal Cartan connection as above.

\begin{defi}
We shall call $(M,H)$ a manifold with a free $n$-distribution. It is always of dimension $m = n(n+1)/2$.
\end{defi}

\subsection{Harmonic curvature} \label{har:cur}

A Cartan connection is said to be \emph{torsion-free} if the curvature $\kappa$ seen as a function $\mc{P} \to \wedge^2 \mf{p}^{\perp} \otimes \mf{g}$ is actually a function $\mc{P} \to \wedge^2 \mf{p}^{\perp} \otimes \mf{p}$. If $\kappa$ is instead seen as a section of $\wedge^2 T^* \otimes \mc{A}$, torsion-freeness implies that is a section of $\wedge^2 T^* \otimes \mc{A}_{(0)}$.

Since $\partial^* \kappa = 0$, $\kappa$ must map into Ker$(\partial^*)$. This has a projection onto the homology component $H_2(\mf{p}^{\perp}, \mf{g}) =$ Im$(\partial^*)/$Ker$(\partial^*)$. The composition of $\kappa$ with this projection is $\kappa_H$, the \emph{harmonic curvature}. The Bianci identity for normal Cartan connections imply that it is a complete obstruction to integrability (\cite{CartEquiv}). Indeed, paper \cite{BGG} demonstrates that $\kappa_H$ determines $\kappa$ entirely.

Now, Kostant's version of the Bott-Borel-Weil theorem \cite{Kostant} allows us to algorithmically calculate $H_2(\mf{p}^{\perp}, \mf{g})$. The $n = 1$ case is trivial, the $n=2$ and $n=3$ cases have $H_2(\mf{p}^{\perp}, \mf{g})$ contained inside
\be
(\mf{g}_{1} \wedge \mf{g}_{2}) \otimes \mf{g}_{0}.
\ee
The harmonic curvature must lie inside this component, which is of homogeneity three. Both the Bianci identity for normal Cartan connections \cite{CartEquiv}, and the construction of the full curvature from the harmonic curvature \cite{BGG} imply that the other component of the curvature must have higher homogeneity. Since the torsion components have maximal homogeneity three (for $(\wedge^2 \mf{g}_2) \otimes \mf{g}_{-1}$), these two geometries are torsion-free.

For $n \geq 4$, the harmonic curvature is contained inside
\be
(\mf{g}_{1} \wedge \mf{g}_{2}) \otimes \mf{g}_{-2},
\ee
and thus these geometries are never torsion-free (unless they are flat).

\subsection{The Tractor bundle} \label{trac:bun}
The standard Tractor bundle $\mc{T}$ is bundle on which we will be doing most of our calculations. It is defined to be the bundle associated to $\mc{P}$ under the standard representation of $\mf{so}(n+1,n)$, restricted to $\mf{p}$:
\be
\mc{P} \times_P \mathbb{R}^{(n+1,n)}.
\ee

This makes $\mc{T}$ into a rank $2n+1$ bundle. A choice of preferred connection $\nabla$ reduces the structure group of $\mc{T}$ to $\mf{gl}(n)$. Under this reduction,
\be
\mc{T} = H \oplus \mathbb{R} \oplus H^*.
\ee
Changing the the choice of $\nabla$ by a one-form $\Upsilon$ changes this splitting as:
\beqa \label{splitting:change}
\left( \begin{array}{c} v \\ \tau \\ X \end{array} \right) \to \left( \begin{array}{c} v + \tau \Upsilon_1 - \{\Upsilon_2, X \} - \frac{1}{2} (\Upsilon_1(X)) \Upsilon_1 \\ \tau -\Upsilon_1(X) \\ X \end{array} \right),
\eeqa
where $\Upsilon_1$ and $\Upsilon_2$ are the components of $\Upsilon$ under the splitting of $T^* = T^*_1 \oplus T^*_2$ given by $\nabla$.

This demonstrates that the inclusions $H^* \subset \mathbb{R} \oplus H^* \subset \mc{T}$ are well defined, as are the projections $\mc{T} \to H \oplus \mathbb{R} \to H$. The projection that we shall be using the most often is $\pi^2: \mc{T} \to H$.

The metric $h$ on $\mc{T}$, of signature $(n+1, 1)$ is given in this splitting by:
\be
h \ \big( \left( \begin{array}{c} v \\ \tau \\ X \end{array} \right) , \left( \begin{array}{c} w \\ \nu \\ Y \end{array} \right) \big) = \frac{1}{2} (w(Y) + v(X) + \tau \nu),
\ee
while the Tractor derivative in the direction of any $Z \in \Gamma(T)$ is $\onab_Z = Z + \nabla_Z + \rP(Z)$, or, more explicitly,
\be \label{tractor:con}
\onab_Z \left( \begin{array}{c} v \\ \tau \\ X \end{array} \right) = \left( \begin{array}{c} \nabla_Z v + \tau \rP(Z)_1 - \{\rP(Z)_2 , X\} \\ \nabla_Z \tau - v(Z_{-1}) - \rP(Z)_1(X) \\ \nabla_Z X +  \tau Z_{-1} + \{ Z_{-2}, v \} \end{array} \right).
\ee

\section{Fefferman constructions}
Consider a parabolic geometry $(M, \mc{P}, \omega)$ derived from the homogeneous model $G/P$. Assume that there is an inclusion $G \subset \widehat{G}$ with a parabolic inclusion $\wh{P} \subset \wh{G}$ such that $\wh{P} \cap G \subset P$. Assume further that the inclusion $G/(\wh{P} \cap G) \subset \wh{G} / \wh{P}$ is open. Then we may do the \emph{Fefferman construction} on this data. See for example \cite{CR} for details of the original construction.

Define $\wh{M}$ as $\mc{P} / (\wh{P} \cap G)$. The inclusion $(\wh{P} \cap G) \subset \wh{P}$ defines a principal bundle inclusion $i: \mc{P} \hookrightarrow \wh{\mc{P}}$ over $\wh{M}$. Since $\mf{g} \subset \wh{\mf{g}}$, we may extend $\omega$ to a section $\omega'$ of $(T\wh{\mc{P}}^* \otimes \wh{\mf{g}})|_{\mc{P}}$ by requiring that $\omega'(\sigma_A) = A$, for any element $A \in \wh{\mf{g}}$ and $\sigma_A$ the fundamental vector field on $\wh{\mc{P}}$ defined by $A$. We may further extend $\omega'$ to all of $\wh{\mc{P}}$ by $\wh{P}$-equivariance.

Since the inclusion $G/(\wh{P} \cap G) \subset \wh{G} / \wh{P}$ is open, the inclusion $\mf{g} \subset \wh{\mf{g}}$ generates a linear isomorphism $\mf{g}/(\wh{\mf{p}} \cap \mf{g}) \to \wh{\mf{g}} / \wh{\mf{p}}$. At any point $u \in \mc{P}$, $\omega$ is a linear isomorphism $T\mc{P}_u \to \mf{g}$. The previous condition ensures that $\omega'$ is a linear isomorphism $T\wh{\mc{P}}_u \to \wh{\mf{g}}$. This condition extends to all of $\mc{P}$, then to all of $\wh{\mc{P}}$ by equivariance. Consequently $\omega'$ is a Cartan connection.

Dividing out by $P / (\wh{P} \cap G)$ makes $\wh{M}$ into a principal bundle over $M$. It is then easy to see that $\omega'$ is invariant along the vertical vectors of $\wh{M}$ and projects to $\omega$ on $M$. Thus $\omega'$ has holonomy group contained in $G$.

\begin{rem}
The Fefferman construction implies nothing about the relative normalities of $\omega$ and $\omega'$.
\end{rem}
\subsection{Almost-spinorial structures} \label{spin}

There is an evident inclusion of $SO(n+1,n)$ into $SO(n+1,n+1)$. In terms of Dynkin diagrams,
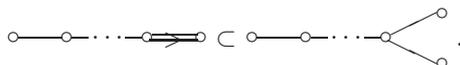
\be
\Balgebra \subset \Dalgebra \ .
\ee

\begin{prop}
There exists a Fefferman construction for this inclusion, where $\wh{M} = M$. In terms of Dynkin diagrams with crossed nodes (see \cite{CartEquiv}), this is
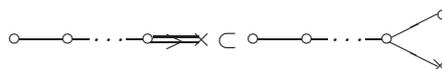
\be
\Bpalgebra \subset \Dpalgebra \,
\ee
and the other parabolic geometry is an almost-spinorial geometry (see \cite{Invar}). This Fefferman construction was discovered by Doubrov and Slov\'ak, \cite{feffsame}.
\end{prop}
\begin{proof}
The homogeneous model for the almost-spinorial geometry is $\wh{G} / \wh{P}$ where $\wh{G} = SO(n+1,n+1)$ and $\wh{P}$ is the stabilizer of an isotropic $n+1$ plane. The homogeneous model for a free $n$-distribution are given by $G = SO(n+1,n)$ and $P$ the stabilizer of an isotropic $n$ plane. Since the space $\mbb{R}^{(n+1,n)} \subset \mbb{R}^{(n+1,n+1)}$ must be transverse to every isotropic $n+1$ plane, $\wh{P} \cap G = P$. The open inclusion for the Fefferman condition is equivalent with the statement that $\mf{g}$ and $\wh{\mf{p}}$ are transverse inside $\wh{\mf{g}}$. A simple comparisons of the ranks of all the algebras involved demonstrates that this is the case. This allows us to do the Fefferman construction.

Since $\wh{P} \cap G = P$,
\be
\wh{M} = \wh{\mc{P}} / \wh{P} = \mc{P} / \wh{P} \cap G = \mc{P} / P = M.
\ee
So this almost spinorial structure is on the same manifold as the free $n$-distribution.
\end{proof}

There is another Fefferman construction that may initially seem relevant here; that given by the inclusion
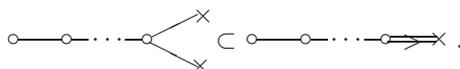
\be
\Dpalgebrajaw \subset \Bpalgebra \ .
\ee
But except when the first algebra is $D_4$ or $D_3$, parabolic geometries of the the first type are flat if regular and normal (since all their harmonic curvatures have zero homogeneity, see Kostant's version of the Bott-Borel-Weil theorem \cite{Kostant}). The case of $D_3$ will be dealt with in Section \ref{CR:Feff} -- it is the CR, standard Fefferman construction.

\section{Free $3$-distributions} \label{4:3}
These are the geometries detailed by Bryant in \cite{bryskew}. They have two properties that distinguish them from the general free $n$-distribution behaviour. First of all, they are torsion free (Section \ref{har:cur}). Secondly, the almost spinorial Fefferman construction of Section \ref{spin} is given by $\trione$. However, triality implies that
\be
\trione  \ \ \cong \ \tritwo \ ,
\ee
i.e.~that the almost-spinorial structure is actually a conformal structure, whenever the $SO(4,3)$ structure lifts to a $Spin(4,3)$ structure. This is always true locally.

Paper \cite{bryskew} details the Fefferman construction explicitly. He further shows that if the Tractor connection for the free $3$-distribution is regular and normal, the conformal Tractor connection must be normal as well (regularity is automatic since the conformal parabolic is $|1|$-graded). The holonomy of that conformal Tractor connection must evidently reduce to $Spin(3,4)$.

In fact, the conformal structure is determined by the filtration of $T$ coming from the Tractor connection of the $3$-distribution (see next section). Consequently this local lift globalises for all free $3$-distributions.

\begin{prop}
Conversely, if the normal conformal holonomy of a six manifold $M$ of split signature reduces to $Spin(4,3) \subset SO(4,4)$, this manifold is locally the Fefferman construction of a regular normal free $3$-distribution.
\end{prop}
\begin{proof}
Recall that the spin representation of $Spin(4,3)$ is $\mbb{R}^(4,4)$, so the above proposition makes sense. Set $\wh{G} = SO(4,4)$, with $\wh{P}$ being $CO(3,3) \rtimes \mbb{R}^{(3,3)}$, the conformal parabolic (defined as the stabiliser of a null-line in $\mbb{R}^{(4,4)}$). As before, $G = Spin (4,3)$ and $P = GL(3) \rtimes \mbb{R}^{3} \rtimes \wedge^2 \mbb{R}^3$. Let $\onab^c$ and $\omega^c$ be the normal conformal Tractor and Cartan connections.

Let $\wh{\mc{P}}$ be the conformal $\wh{P}$ bundle, $\wh{\mc{P}} \subset \wh{\mc{G}}$ with $\wh{\mc{G}}$ the full structure bundle for $\omega^c$. The holonomy reduction implies that there exists a $G$-bundle $\mc{G} \subset \wh{\mc{G}}$ such that $\omega^c$ reduces to a principal connection on $\mc{G}$.

The action of $Spin(4,3)$ on the null-lines of $\mbb{R}^{(4,4)}$ is transitive under the spin representation, consequently there is only one type of intersection between $G$ and $\wh{P}$ inside $\wh{G}$. By \cite{bryskew}, this is $G \cap \wh{P} = P$. This means that $\mc{G} \cap \wh{\mc{P}} = \mc{P}$, a $P$-bundle, so $\omega^c$ reduces to a free $3$-distribution Cartan connection -- call it $\omega$.

It remains to show that this Cartan connection is normal. Looking at the homogeneous model, the conformal structure comes from the fact there is a unique conformal class of $\wh{P}$-invariant inner products on $\wh{\mf{g}}/\wh{\mf{p}}$. This implies there is a unique conformal class of $P$ invariant inner products on $\mf{g}/\mf{p}$.

Since $TM = \mc{G} \times_P (\mf{g}/\mf{p})$, this means that the conformal structure on $TM$ depends only on the component of $\omega$ mapping into negative homogeneity -- this is the soldering form, $\omega_-$ (see next section for the geometric details of this).

The curvature $\kappa^c$ of $\omega^c$ can be seen as a $\wh{P}$-invariant map from $\wh{\mc{P}}$ to $\wedge^2 (\wh{\mf{g}}/\wh{\mf{p}})^* \otimes \wh{\mf{g}}$. Similarly, the curvature $\kappa$ of $\omega$ is a $P$-invariant map from $\mc{P}$ to $\wedge^2 (\mf{g}/\mf{p})^* \otimes {\mf{g}}$. On $\mc{P}$, these two curvatures are related by the commuting diagram:
\be
\begin{array}{ccc}
\wedge^2 \mf{g}/\mf{p} &\stackrel{\kappa}{\longrightarrow}& \mf{g} \\
\uparrow & & \downarrow \\
\wedge^2 \wh{\mf{g}}/\wh{\mf{p}} &\stackrel{\kappa^c}{\longrightarrow}& \wh{\mf{g}}.
\end{array}
\ee
Since $\omega^c$ is normal, it is torsion free (see \cite{TCPG} or \cite{mecon}), implying that it maps into $\wh{\mf{p}}$. This means that $\omega$ also maps into $\mf{p}$ -- so is also torsion-free. This means that $\kappa$ is of homogeneity $\geq 2$, consequently -- since $\partial^*$ respects homogeneity -- $\partial^* \kappa$ is of homogeneity $\geq 2$.

Now, by \cite{capslo}, any Cartan connection $\omega$ with curvature $\kappa$ such that $\partial^* \kappa$ is of homogeneity $\geq l \geq 0$ differs from the normal Cartan connection $\omega'$ by a section $\Phi \in \Omega^1(\mc{P},\mf{g})$ of homogeneity $\geq l$.

So here we have $\omega + \Phi = \omega'$, with $\omega'$ normal and $\Phi$ of homogeneity $\geq 2$. This means that $\omega'$ and $\omega$ have the same soldering form (as the soldering form is of homogeneity $\leq 1$, as the algebra is two-graded), thus that the conformal structure that they both generate are the same. Since the conformal Fefferman construction for $(\mc{P}, \omega')$ must be normal (since $\omega'$ is), it must \emph{be} $(\wh{\mc{P}}, \omega^c)$. This means that $\Phi = 0$, hence that $\omega$ is normal.
\end{proof}

\subsection{Geometric equivalence}

Given a free $3$-distribution on manifold $M$, the conformal structure can be recovered directly from the decomposition of $T \cong T_{-2} \oplus H$ given by any Weyl structure. Let $\sigma$ be any local never-zero section of $\wedge^3 H$. Then there is a map $g: T_{-2} \otimes H \to \wedge^3 H$ given by the isomorphism $T_{-2} \cong \wedge^2 H$. Extend $g$ to a section of $(\odot^2 T^*) \otimes \wedge^3 H$ by the inclusion $T_{2}^* \otimes H^* \subset \odot^2 T^*$. Then $g\sigma^{-1}$ is a metric on $M$. This depends on the choice of the section $\sigma$, so actually defines a conformal structure. It is then easy to see that $g$ is invariant under the action of a one-form $\Upsilon$, (as $g(U+Y,X) = g(U,X) + g(X,Y) = g(U,X)$ for any sections $X$ and $Y$ of $H$ and any section $U$ of $T_{-2}$). So this conformal structure does not depend on the choice of Weyl structure, only on the filtration of $T$ (which depends on the Cartan connection).

It can be instructive to construct the conformal structure $[g]$ directly from the distribution, without having to construct the full Cartan connection. As seen above, it suffices to construct a single compatible transverse distribution $T_{-2}$ -- thus to find a single preferred connection $\nabla$. Before doing so, we need to note the relationship between the torsion of $\nabla$, $\rP$ and the curvature $\kappa$ of $\onab$. For sections $X$ and $Y$ of $T$, this relation is:
\be
\kappa_{-}(X,Y) = Tor^{\nabla}(X,Y) + \{X,Y\} + \{\rP(X),Y\}_{-} - \{\rP(Y), X\}_{-}.
\ee
See paper \cite{CapWeyl} for more details. Since our Tractor connection is torsion-free, the left hand term vanishes. We will be looking specifically at the homogeneity $(1,1,-1)$, $(1,1,-2)$ and $(1,2,-2)$ components of $Tor^{\nabla}$. In all these homogeneities, $\{\rP(X),Y\}_{-} - \{\rP(Y), X\}_{-}$ vanishes, giving us:
\be
Tor^{\nabla}(X,Y) = -\{X,Y\}.
\ee

\begin{lemm}
The transverse distribution for a preferred connection $\nabla$ for a generic $3$-distribution $H$ depends only on the torsion-freeness of $\onab$ and the values of $\nabla_A B$, for all $A$ and $B$ sections of $H$.
\end{lemm}
\begin{proof}
Since $\nabla$ is torsion-free, the $H$ component of $[A,B]$ is $\nabla_A B - \nabla_B A$ (homogeneity $(1,1,-1)$). Thus $[A,B] -\nabla_A B + \nabla_B A$ is a section of $T_{-2}$, thus defining this bundle.
\end{proof}

\begin{prop} \label{partial:prop}
Let $\nabla$ be a partial connection, differentiating sections of $H$ along directions in $H$. It may be extended, via the algebra{\"{\i}}c bracket $\{,\}$, to differentiate sections of $T/H$ along directions in $H$. Let $q$ be the projection $T \to T/H$. Define a transverse distribution $T_ {-2}$ with a projection $\Pi: T \to H$. Take $\nabla$ to be torsion-free, in the sense that for $A$  and $B$ sections of $H$ and $Z$ any section of $T$,
\beqa
\label{conformal:one} \nabla_A B - \nabla_B A - \Pi([A,B]) &=& 0,\\
\label{conformal:two} \nabla_A q(Z) - q([A,Z]) + \{A,\Pi(Z)\} &=& 0,
\eeqa
(the top term is of homogeneity $(1,1,-1)$; the bottom one is a mix of $(1,1,-2)$ and $(1,2,-2)$). Then $\nabla$ and $\Pi$ (and hence $T_{-2}$) are entirely determined by the action of $\nabla$ on the line-bundle $\wedge^3 H^*$.
\end{prop}
\begin{proof}
Let $A, B, C$ and $D$ be sections of $H$. Inserting $Z = [B,C]$ into equation (\ref{conformal:two}) gives
\beqa \label{p:con}
\{\nabla_A B, C\} + \{B, \nabla_A C\} - q([A,[B,C]]) + \{A, \nabla_B C - \nabla_C B\} =0,
\eeqa
using (\ref{conformal:one}). Let $\sigma$ be a never-zero section of $\wedge^3 H^*$; since the bracket gives an isomorphism $T/H \cong \wedge^2 H$, $\sigma$ is also a section of $H^* \otimes (T/H)^*$, a non-degenerate pairing between $H$ and $T/H$. The exact relationship between the two forms of $\sigma$ is:
\be
\sigma(A,B,C) = \sigma(A, \{B,C\}).
\ee
We shall shift freely between the two definitions.

We wish to estimate $R = \sigma(\nabla_A B,C,D) = \sigma(\nabla_A B, \{C,D\})$, and through it, $\nabla_A B$. First,
\beqa \label{A:eqa}
\nonumber A\cdot\sigma(B,C,D) &=& (\nabla_A \sigma)(B,C,D) + \sigma(\nabla_A B,C,D) + \sigma(B,\nabla_A C,D) + \sigma(B,C,\nabla_A D) \\
&=& (\nabla_A \sigma)(B,C,D) + R + \sigma(B, \nabla_A \{C,D\}) \\
\nonumber &=& (\nabla_A \sigma)(B,C,D) + R + \sigma(B, q([A,[C,D]])) - \sigma(B,A,\nabla_C D - \nabla_D C),
\eeqa
using equation (\ref{p:con}). All of these terms can be computed from the data, apart from $\sigma(B,A,\nabla_C D - \nabla_D C)$. We next need to note that
\beqa \label{CD:eqa}
\nonumber -C \cdot \sigma(D,B,A) + D\cdot\sigma(C,B,A) &=& -(\nabla_C \sigma)(D,B,A) + (\nabla_D \sigma)(C,B,A) \\
&& -\sigma(D, \nabla_C \{B,A\}) + \sigma(C, \nabla_D \{B,A\}) \\
\nonumber && + \sigma(A,B,\nabla_C D - \nabla_D C).
\eeqa
While on the other hand:
\beqa \label{B:eqa}
\nonumber 2 B \cdot \sigma(C,D,A) &=& 2 (\nabla_B \sigma)(C,D,A) + 2 \sigma(A, \nabla_B \{C,D\})  - \sigma(D,C,\nabla_A B) + \sigma(C,D,\nabla_A B) \\
\nonumber && - \sigma(D,C, \nabla_B A - \nabla_A B) + \sigma(C,D, \nabla_B A - \nabla_A B) \\
\nonumber &=& 2 (\nabla_B \sigma)(C,D,A) + 2R + 2 \sigma(A, q([B,[C,D]])) - 2\sigma(A,B, \nabla_C D - \nabla_D C) \\
&& - \sigma(D, q([C,[B,A]])) + \sigma(D,\nabla_C \{B,A\}) \\
\nonumber &&+ \sigma(C, q([D,[B,A]])) - \sigma(C,\nabla_D \{B,A\})
\eeqa
Adding equations (\ref{A:eqa}), (\ref{CD:eqa}) and (\ref{B:eqa}), and rearranging gives the Levi-Civita-like formula:
\be
3 \sigma(\nabla_A B, C, D) &=& A\cdot \sigma(B,C,D) + 2 B\cdot\sigma(A,C,D) + C\cdot\sigma(A,B,D) -D\cdot\sigma(A,B,C) \\
&& -(\nabla_A \sigma)(B,C,D) - 2 (\nabla_B \sigma)(A,C,D) - (\nabla_C \sigma)(A,B,D) + (\nabla_D \sigma)(A,B,C) \\
&& +\sigma(B,q([A,[D,C]])) -2\sigma(A,q([B,[C,D]])) \\
&& + \sigma(D,q([C,[B,A]])) - \sigma(C,q([D,[B,A]])).
\ee
Since $\sigma(\nabla_A B,C,D)$ is entirely determined by the data, the non-degeneracy of $\sigma$ implies that $\nabla_A B$ is also thus determined. Then, by equation (\ref{conformal:one}), $\Pi$ is also thus determined and consequently, so is $T_{-2}$.
\end{proof}

\begin{cor}
Any partial connection $\nabla$ obeying the properties of the previous proposition can be extended into a preferred connection for the $3$-distribution. This extention is unique, up to the action of $T_2^*$. Consequently, any such partial connections define a transverse $T_{-2}$ compatible with the Cartan connection and hence define the conformal structure for that $3$-distribution.
\end{cor}
\begin{proof}
Any partial connection differentiating $L  = \wedge^3 H^*$ in the directions tangent to $H$ may be extended to a full connection differentiating $L$ along all of $T$. The space of full connections extending a given partial connection is parametrised by $T^*_2$.

The preferred connections for the Cartan connection are an affine space, modelled on $T^*$. Full connections on $L$ also form an affine space, similarly modelled on $T^*$. One can easily check, using the change of connection formula in equation (\ref{change:con}), that the map from preferred connections to connections on $L$ is injective (hence bijective).

Let $\nabla$ be a partial connection on $H$ in the sense of Proposition \ref{partial:prop}. Let $\nabla'$ be a preferred connection, such that its action on $L$ extends that of $\nabla$. Since all preferred connections are torsion-free in the sense of Proposition \ref{partial:prop}, the values of the derivatives of $H$ along $H$ via $\nabla'$ are entirely determined by the derivative of $L$ along $H$ via $\nabla'$. Consequently, for any sections $A$ and $B$ of $H$,
\be
\nabla_A B = \nabla'_A B.
\ee
\end{proof}

The converse procedure (constructing the $3$-distribution from a conformal structure with conformal Tractor holonomy algebra $\mf{spin}(4,3)$) can more easily be seen from algebra{\"i}c considerations. The algebra $\mf{spin}(4,3) \subset \mf{so}(4,4)$ is defined as preserving a generic four-form $\lambda$ on $V = \mbb{R}^{(4,4)}$, see \cite{baumspin}. Let $\mc{T}^C$ be the standard conformal Tractor bundle on $M$ (see \cite{TBIPG} or \cite{mecon} for more details on conformal geometries). If the conformal Tractor connection $\onab^C$ has holonomy algebra reducing to $\mf{spin}(4,3)$, then there exists a generic four-form $\nu \in \Gamma(\wedge^4 \mc{T}^C)$ such that
\be
\onab^C \nu = 0.
\ee
There is a natural projection on $\mc{T}_C$, coming from its filtration
\be
\mc{T}_C  \to \mc{E}[1] \ltimes T[-1] \to \mc{E}[1].
\ee
Here $\mc{E}[1]$ is a density bundle, $\mc{E}[\alpha] = (\wedge^6 T^*)^{\frac{\alpha}{-6}}$, and $T[-1] =T \otimes \mc{E}[-1]$. This implies that there is a well defined projection $\pi: \wedge^4{\mc{T}_C} \to (\wedge^3 T)[-2]$. It turns out that $\pi(\nu)$ is decomposable, and so defines a distribution $H^*$ of rank three in $T^*[2/3]$. Since a distribution is unchanged by a change of scale, this is actually a distribution in $T^*$, with dual distribution $H \subset T$. This $H$ is precisely that defining the Bryant structure; the maximal non-integrability derives from the properties of $\nu$ and $\onab^C$.

\subsection{$G_2'$ structures}
The most natural subgroup of $Spin(4,3) \subset SO(4,4)$ is $G_2'$.
\begin{defi}
The group $G_2'$ has several equivalent definitions:
\begin{itemize}
\item $G_2'$ is the subgroup of $Spin(3,4)$ that preserves a given non-isotropic element $e$ in $\mbb{R}^{(4,4)}$.
\item $G_2'$ is the subgroup of $SO(3,4)$ that preserves a generic three-form $\theta$ in $\wedge^3 \mbb{R}^{(3,4)}$.
\item $G_2'$ is the automorphism groups of the split Octonions.
\end{itemize}
\end{defi}
We can see the equivalence between these definitions. As an subgroup of $Spin(3,4)$, $G_2'$ must preserve a generic four-form $\lambda$ in $\mbb{R}^{(4,4)}$. Since it also preserves $e$, it must preserve $e^{\perp}$, giving the inclusion $G_2' \subset SO(4,3)$. And it must preserve the generic three-form $\theta = e\llcorner \lambda$.

The split Octonions are eight-dimensional, and let $V$ be the seven dimensional space $V = Im \mbb{O}'$. $\mbb{O}'$ carries a natural inner product $N$, generated from the norm $N(x,x) = x \overline{x}$. This quadratic form is multiplicative, and is of signature $(4,4)$. The identity element $1 \in \mbb{O}'$ is of positive norm squared, and is orthogonal to $V$; thus $V$ is of signature $(3,4)$. Since automorphisms of the split Octonions must preserve real multiplication, $G_2¸$ preserves $V = 1^{\perp}$. To get the three form $\theta$, we need to use the fact that $\mbb{O}'$ is alternative; meaning that the alternator
\be
[x,y,z] = (xy)z - x(yz)
\ee
is totally anti-symmetric in its three entries. We can use $N$ to make $[,,]$ into an element of $\wedge^3 V^* \otimes V^*$; it turns out to be skew in all four entries. Then $\theta = * [,,]$, where $*$ is the Hodge star on $V$ generated by $N$. Alternatively, it can be defined directly as:
\be
\theta(x,y,z) = N(xy,z),
\ee
skew in all three arguments by the properties of $\mbb{O}'$. In this format, it is obvious that $\theta$ and $N$ allow one to reconstruct the split Octonionic multiplication, giving the equivalence between these three definitions. We shall be using the properties of the split Octonions in the rest of this section.

\subsubsection{Dual distributions}
Now assume that our free $3$-distribution has a normal Tractor connection $\onab$ with a holonomy reduction to $G_2'$. By the conformal Fefferman construction, the conformal structure will be given by a manifold that is conformally Einstein and whose metric cone carries a $G_2'$ structure (see \cite{meein}). Such manifolds do exist -- for instance, $SL(3, \mbb{R}) /T^2$ where $T^2$ is a maximal torus, is one example \cite{brymetric}. Here, the free $3$-distribution would be chosen at $Id \in SL(3, \mbb{R})$ as the span of
\be
H_{Id} = \left( \begin{array}{ccc} 0 & a &0 \\ 0&0&b \\ c&0&0 \end{array} \right),
\ee
and extended to the whole manifold by Lie multiplication. Note that $\{ H_{Id}, \mf{T}^2 \} \subset H_{Id}$, for $\mf{T}^2$ the tangent space to the maximal torus, so this extension is well defined.

It is not unique, however. We could have used the transpose $H_{Id}^t$ instead. Note that $H_{Id}^t = \{H_{Id}, H_{Id} \}$, while $H_{Id} = \{H_{Id}^t, H_{Id}^t\}$. This will be an important property for all $G_2'$ structures on a free $3$-distribution.

\begin{prop}
There are two orbits of isotropic $3$-planes in $\mbb{R}^{(4,3)}$ under the action of $G_2'$. -- one open, one closed. The closed orbit is distinguished by the fact that $\theta(x,y,z) = 0$ for all elements in an isotropic $3$-plane inside this orbit.
\end{prop}
\begin{proof}
Let $B \subset Im \mbb{O}'$ be an isotropic $3$-plane. For any element $x$ of $B$, $0 = N(x,x) = x \times \overline{x} = - x \times x$. Isotropy implies that for any elements $x$ and $y$ of $B$,
\be
0 &=& 2 N(x,y) = N(x+y,x+y) - N(x,x) - N(y,y) \\
&=& (x+y) \times (\overline{x} + \overline{y}) - x \times \overline{x} - y \times \overline{y} \\
&=& -(x+y) \times (x+y) \\
&=& -x \times y - y \times x.
\ee
Thus an isotropic $3$-plane is a subset of $Im \mbb{O}'$ where all elements anti-commute and square to zero. Since the split Octonions are alternative, the multiplicative span of any two elements is associative. Hence
\be
(xy)(xy) = x(yx)y = -x(xy)y = -(xx)(yy) = 0.
\ee
This is true for any elements $x$, $y$ in B. Thus $C = B \times_{ \mbb{O}'} B$ is isotropic. Because the elements of $B$ anti-commute, $C$ is also pure imaginary. There are two situations to be considered:
\begin{enumerate}
\item There exists a basis $\{x,y,z\}$ for $B$ such that $\theta(x,y,z) = \frac{1}{2}$.
\end{enumerate}
The set of all such $B$ is open and dense in the set of all isotropic $3$-planes. We aim to show $G_2'$ is transitive on this set.

\begin{lemm}
The elements span of $x$, $y$ and $z$ under split Octonionic multiplication generate all of $\mbb{O}'$.
\end{lemm}
\begin{lproof}
Here the relation
\beqa \label{77}
\frac{1}{2} = \theta(x,y,z) = N(xy,z) = \frac{1}{2}((xy)z + z(yx)).
\eeqa
implies that $xy$ is orthogonal to $x$, $y$, but not to $z$. We may cyclically permute $x$, $y$ and $z$ here, thus demonstrating that $xy$, $yz$ and $zx$ are linearly independent, hence that $C$ is of dimension three. Define $a$ as
\be
a = (xy)z - z(xy) = [xy,z].
\ee
Evidently, $\overline{a} = -a$, so $a \in Im \mbb{O}'$. Set $\mc{F} = \{x,y,z,xy,yz, a\}$. Using the fact that the multiplicative span of two elements is associative and the fact that $\theta$ is skew, we can see that all the elements of $\mc{F}$ are orthogonal to each other, with the exception of $N(a,a)$ and
\be
N(xy,z) = N(yz,x) =  N(zx,y) =  \frac{1}{2}.
\ee
Thus $(xy+z)$ and $(xy-z)$ are of norm-squared one and minus one, respectively. Since the split-Octonions are a division algebra, $a = (xy+z)(xy-z)$ must be of norm-squared minus one, so
\be
-1 = N(a,a) = -a^2.
\ee
Similarly, $(xy + z) \times (yz + x)) = (xy)(yz) + zx$ must be of norm one. Both elements $zx$ and $(xy)(yz)$ are isotropic, as products of elements on isotropic $3$-planes; thus
\be
\frac{1}{2} = N((xy)(yz),zx) = \theta(xy,yz,zx).
\ee
Now $(xy)(yz)$ is imaginary, and orthogonal to $a$ (as its product with $(xy)$ vanishes, so $N((xy)(yz),(xy)z) - N((xy)(yz),z(xy)) = 0-0$). The above relation thus implies that $(xy)(yz) = y$. Similarly, $(yz)(zx) = z$ and $(zx)(xy) = x$.

There is another expression for $\theta$; it is given as half the commutator $\frac{1}{2} [,]$, combined $N$. Thus we may say that:
\be
\theta(xy,z,a) = N(\frac{1}{2}[xy,z],a) = \frac{1}{2}N(a,a) = -\frac{1}{2}.
\ee
Now $N(za,y) = N(yz,a) = 0$, and similarly $N(za,z) = N(za,zx) = N(za,yz) = 0$. Furthermore, $za = z(xy)z$, which is clearly imaginary. We may therefore conclude that $za= -z$, and, similarly:
\be
az = -za = z \ \ \textrm{and} \ \ (xy)a = -a(xy) = xy.
\ee
Now let $\{x', y', z'\}$ be another basis for $B$ with $\theta(x',y',z') = \frac{1}{2}$. Then by similar reasoning to above, $a' = [x'y',z']$ is orthogonal to $B \oplus C$, and is pure imaginary and of norm squared minus one. Thus $a' = \pm a$. In other words, there is a continuous function from the set of frames of $B$ with this scale, to the set $\{a,-a\}$. Since the first set is connected, this continuous function must be constant. Therefore $a = [yz,x] = [zx,y]$. This allows us to calculate the last remaining multiplicative terms, and gives the full multiplication table for $\mc{F}$ (see table \ref{table:one}). Hence all of split-Octonionic multiplication can be derived from $x$, $y$ and $z$.

\begin{table}[htbp]
\begin{center}
\begin{tabular}{|c||c|c|c|c|c|c|c|}
\hline
$\times$& $a$&$x$&$y$&$z$&$yz$&$zx$&$xy$ \\
\hline
\hline
$a$&$-1$&$x$&$y$&$z$&$-yz$&$-zx$&$-xy$ \\
\hline
$x$&$-x$&$0$&$xy$&$-zx$&$\frac{1}{2}(1-a)$&$0$&$0$ \\
\hline
$y$&$-y$&$-xy$&$0$&$yz$&$0$&$\frac{1}{2}(1-a)$&$0$ \\
\hline
$z$&$-z$&$zx$&$-yz$&$0$&$0$&$0$&$\frac{1}{2}(1-a)$ \\
\hline
$yz$&$yz$&$\frac{1}{2}(1+a)$&$0$&$0$&$0$&$z$&$-y$ \\
\hline
$zx$&$zx$&$0$&$\frac{1}{2}(1+a)$&$0$&$-z$&$0$&$x$ \\
\hline
$xy$&$xy$&$0$&$0$&$\frac{1}{2}(1+a)$&$y$&$-x$&$0$ \\
\hline
\end{tabular}
\end{center}
\caption{Split-Octonionic multiplication for $\mc{F}$}\label{table:one}
\end{table}
\end{lproof}

So call $\{x,y,z\}$ a \emph{split Octonionic triple}. Any element $g$ of $G_2'$ is entirely determined by $\{g(x),g(y),g(z)\}$. Conversely, for any two split Octonionic triple, the set map sending one triple to the other extends to an automorphism of $Im(\mbb{O}')$ respecting split Octonionic multiplication; by definition, this is an element of $G_2'$.

Moreover, if $G_{B} \subset SL(7,\mbb{R})$ is the stabiliser of $B$, $G_2' \cap G_B$ is the permutation group of the Octonionic triples in $B$ -- consequently $G_2' \cap G_B = SL(3,\mbb{R})$, since $\theta$ is a volume form on $B$.

\begin{enumerate}
\setcounter{enumi}{1}
\item For all $x,y,z \in B$, $\theta(x,y,z) = 0$.
\end{enumerate}
The set of all such $B$ is closed in the set of all isotropic $3$-planes, complementary to the previous orbit, and with empty interior. We aim to show $G_2'$ is transitive on this set.

Though all its properties can be deduced from $\theta$ and $*\theta$, it often helps to work with an explicit description of split Octonion multiplication. Here is one, due to Zorn. A split Octonion is represented by the ``matrix''
\be
x= \left( \begin{array}{cc} a & \bf{v} \\ \bf{w} & b \end{array}\right)
\ee
with $a$ and $b$ real numbers and $\bf{v}$ and $\bf{w}$ vectors in $\mbb{R}^3$. The norm squared $N(x,x)$ is the ``determinant'' $ab - \bf{v} \cdot \bf{w}$. Multiplication is given by
\be
\left( \begin{array}{cc} a & \bf{v} \\ \bf{w} & b \end{array}\right) \times \left( \begin{array}{cc} a' & \bf{v}' \\ \bf{w}' & b' \end{array}\right) = \left( \begin{array}{cc} aa' + \bf{v}\cdot\bf{w}' & a\bf{v}' + b'\bf{v} + \bf{w} \wedge \bf{w}' \\ a'\bf{w} + a \bf{w}' - \bf{v} \wedge \bf{v}' & bb' +\bf{v}' \cdot \bf{w} \end{array}\right).
\ee
With $\cdot$ and $\wedge$ the ordinary dot and cross products on $\mbb{R}^3$. The imaginary split Octonions are those where $a  = -b$.

Now, as in case when $\theta$ did not degenerate on $B$, $(xy)(xy) = 0$, and $C = B \times_{\mbb{O}'} B$ is isotropic. However
\be
0 = \theta(x,y,z) = N(xy,z),
\ee
implying that $C \subset B^{\perp}$. Since $C$ is isotropic, $C \subset B$. An inspection of the explicit form of split Octonion multiplication demonstrates that there does not exist a three plane on which $\times_{\mbb{O}'}$ is totally degenerate. So $C \neq 0$. Pick a $z \neq 0$ in $C$. Now $z = xy$ for elements $x$ and $y$ in $B$. Since elements of $B$ square to zero, $x \neq y$. Since the multiplicative span of any two elements is associative, $z \neq 0$ and $zy = zx = 0$, we know that $z \neq x$ and $z\neq y$. Furthermore, $z$ cannot be in the linear additive span of $x$ and $y$, since $x (r_1 x + r_2 z) = r_1 xx + r_2 x(xy) = 0$ for all real $r_j$. So $x$, $y$ and $z$ form a basis for $B$, and the relations
\be
xy = z, yx = -z, xz =0, yz = 0, xx = yy = zz = 0,
\ee
determine multiplication on $B$ completely. In fact, $B$ is determined by $z$. This can be seen from the fact that $G_2'$ is transitive on the set of isotropic element of $Im \mbb{O}'$, so we may set
\be
z = \left( \begin{array}{cc} 0 & e_1 \\ 0 & 0 \end{array} \right),
\ee
where $e_1, e_2, e_3$ is a basis for $\mbb{R}^3$. Then the two sided kernel of the multiplications $\times z, z \times : Im \mbb{O}' \to \mbb{O}'$ is spanned by
\be
z, \left( \begin{array}{cc} 0 & 0 \\ e_2 & 0 \end{array} \right), \left( \begin{array}{cc} 0 & 0 \\ e_3 & 0 \end{array} \right).
\ee
Since $B$ is in the two-sided kernel of multiplication by $z$, and is isotropic, it must be precisely the span of these elements. Since $B$ is determined by $z$, and since $G_2'$ is transitive on isotropic elements of $Im \mbb{O}'$, $G_2'$ must be transitive on the set of isotropic $3$-planes $B$ on which $\theta$ vanishes.

By the above, the subgroup of $G_2'$ that preserves $B$ is the same subgroup that preserves $z$: in other words, it is the nine-dimensional algebra (see \cite{kat}):
\beqa \label{Poct}
\mf{gl}(2) \oplus \mbb{R}^2 \oplus (\wedge^2 \mbb{R}^2) \oplus (\wedge^2 \mbb{R}^2) \otimes \mbb{R}^2
\eeqa
\end{proof}

\begin{prop}
Let $M$ be a free $3$-distribution manifold with normal Tractor connection $\onab$, with the holonomy group of $\onab$ reducing to $G_2'$, with corresponding principal bundle $\mc{G}_2' \subset \mc{G}$. Define
\be
\mc{B} = \mc{G}'_2 \times_{G_2'} \mf{g}'_2 \subset \mc{A}.
\ee
Then on an open, dense subset of $M$, $\mc{B} \cap \mc{A}_{(0)}$ is of rank eight or less. This implies that on this subset, $\mc{B} \cap \mc{A}_{(0)}$ must be an algebra bundle modeled on the algebra $\mf{sl}(3)$, and that $\B \cap \A_{(1)} = 0$.
\end{prop}
\begin{proof}
We shall demonstrate that the projection $\mc{B} \to T$ is surjective on an open dense set of $M$; this then proves the result, as $\mc{B}$ is of rank $14$ and $T$ is of rank $6$, while $\mc{A}_{(0)}$ is the kernel of the projection.

Pick any Weyl structure $\nabla$ to get a splitting of $\mc{A}$. We shall first show that whenever $\mc{B} \cap \mc{A}_{(2)} \neq 0$ (an open condition) the projection onto $T$ is surjective. To see that, pick any nowhere-zero local section $\nu$ of $\mc{B} \cap \mc{A}_{(2)} = \mc{B} \cap T^*_2$, and differentiate it repeatedly. Call the span of these derivatives $\nu^{\onab}$. By definition, for $X_j$ sections of $H = T_{-1}$, $(\onab_{X_1} \ldots \onab_{X_p} \nu)_{p-2} = \{X_1, \{ \ldots, \{ X_p, \nu \} \ldots \}$. This implies that $\nu^{\onab}_1$ is of rank at least two, $\nu^{\onab}_0$ is of rank at least six, and $\nu^{\onab}_{-1}$ and $\nu^{\onab}_{-2}$ are both of rank three -- thus the projection onto $T$ is surjective, as $\onab$ preserves $\mc{B}$. Similar reasoning demonstrates the same result whenever $\mc{B} \cap \mc{A}_{(1)} \neq 0$.

So now assume that $\mc{B} \cap \mc{A}_{(1)} = 0$, thus that the projection $\mc{B} \to T \oplus \mc{A}_0$ is surjective. This means that the projection of $\mc{B}$ onto $\mc{A}_0$ is of rank eight or nine. If it is of rank eight, we are done; if is of rank nine on any open set, there exist a local section $\phi$ of $\mc{B}\cap \mc{A}_{(0)}$ such that $\phi$ projects to the grading section on $\mc{A}_0$. Thus
\be
(\onab_{Z_{-j}} \phi)_{-j} = \{Z_{-j},\phi_0\} = j Z_{-j},
\ee
demonstrating that the map $\mc{B}$ to $T$ is surjective, and hence that $\mc{B} \cap \mc{A}_{(0)}$ is actually of rank eight, on an open, dense subset.

Finally, the fact that $\mc{B} \cap \mc{A}_{(0)}$ must be modeled on the algebra $\mf{sl}(3)$ comes directly from the previous calculations: there are only two ways for $\mf{g}_2'$ and $\mf{p}$ to intersect, either as $\mf{sl}(3)$ (of dimension eight) or as the nine dimensional algebra defined in (\ref{Poct}). Passing to a bundle, the result follows.
\end{proof}

\begin{theo}
Let $M$ be a free $3$-distribution manifold with normal Tractor connection $\onab$, with the holonomy group of $\onab$ reducing to $G_2'$. Then, on an open, dense set of $M$, there is a unique Weyl structure $\nabla$ defined by this information. This Weyl structure determines a splitting of $T = T_{-2} \oplus H$.

Then $H'=T_{-2}$ is a free $3$-distribution, and this information also determines, up to isomorphism, a unique commutative bundle inclusion graph:
\be
\begin{array}{ccccc}
&&\mc{G}&& \\
&\nearrow &&\nwarrow& \\
\mc{P} &&&& \mc{P}' \\
&\nwarrow&&\nearrow& \\
&&\mc{G}_0&&
\end{array},
\ee
where $\mc{P}'$ is a principal $P'$-bundle with $P' \cong P$. Pulling back the Tractor connection $\omega$ from $\mc{G}$ to $\mc{P}'$, we get a Cartan connection on $\mc{P}$, which is the normal Cartan connection for the free $3$-distribution $H'$.

We may iterate this process, since $\omega$ still has holonomy contained in $G_2'$. If we do so, we get back to where we started, with $(H')' = H$. Thus $HA$ and $H'$ can be seen as \emph{dual distributions}.
\end{theo}
\begin{proof}
Since $\onab$ has holonomy contained in $G_2'$, for the rest of the proof we will restrict to the open dense set of $M$ where $\mc{B}_0 = \mc{B} \cap \mc{A}_{(0)}$ is of rank eight. Hence $\mc{B}_0$ is an algebra bundle with the structure of $\mf{sl}(3)$ on each fiber.

Split-Octonionic multiplication gives a well defined subbundle of $\mc{T}$, $K = H^* \times H^*$. Since $K$ and $H^*$ are transverse, the projection $\pi^2$ maps $K$ isomorphically to $H$. Now, there is a unique Weyl structure such that $\mc{T}$ splits as $K \oplus \mbb{R} \oplus H^* $ (basic algebraic manipulation; see \cite{meskew} for a detailed look at this). Then let $\nabla$ be the preferred connection equivalent to this Weyl structure.

Pick a point $u$ of $\mc{P}$. This gives an identification $i_u: \mc{T}_{\pi(u)} \to \mbb{R}^{(4,3)}$. Define $V \subset \mbb{R}^{(4,3)}$ to be the canonical subspace corresponding to $H^* \subset \mc{T}$, and $W \subset \mbb{R}^{(4,3)}$ the space corresponding to $K$. We then define $\mc{G}_0$ as the subbundle of $\mc{P}$ such that
\be
i_v(K_v) = W.
\ee
It is easy to see that this bundle is $G_0$-bundle. Further define $\mc{P}'$ as the subbundle of the full principal bundle $\mc{G}$ on which the preceding property holds. Since $K$ is isotropic, $\mc{P}$ must be a $P'$ bundle, for $P'$ a Lie group conjugate to $P$ with $\mf{p} = \mf{g}_0 + \mf{g}_{-1} + \mf{g}_{-2}$ (we have used the grading section defined by this splitting to split $\mf{g}$). Furthermore, since $\mc{P} \subset \mc{G}$ is defined as the subbundle for which $i_v(H^*_v) = V$, we can see that $\mc{P} \cap \mc{P}' = \mc{G}_0$ and that we have bundle inclusions:
\be
\begin{array}{ccccc}
&&\mc{G}&& \\
&\nearrow &&\nwarrow& \\
\mc{P} &&&& \mc{P}' \\
&\nwarrow&&\nearrow& \\
&&\mc{G}_0&&
\end{array}
\ee
These bundle inclusions depend on the choice only of a $u$, making them invariant up to the action of $P$ -- hence up to isomorphism.

If $\omega$ is the principal connection corresponding to $\onab$, then $\omega$ pulls back to a section of $T\mc{P}' \otimes \mf{g}$, that is $P'$-equivariant and maps vectors generated by $A \in \mf{p}'$ to $A$. Our contention is that this pull-back is a normal Cartan connection.

To do so, first pull back $\omega$ to $\mc{G}_0$ (this pull back factors through pull backs via both $\mc{P}$ and $\mc{P}'$). There, $\omega$ must decompose as
\be
\omega_{-2} + \omega_{-1} + \omega_0 + \omega_1 + \omega_2.
\ee
The $\omega_0$ piece is the component corresponding to $\nabla$ (since $\mc{G}_0$ preserves the splitting of $\mc{T}$). The $\omega_{-2} + \omega_{-1}$ corresponds to the soldering form for $(\mc{P},\omega)$. Finally, $\omega_1 + \omega_2$ is the ``soldering form'' for $(\mc{P}',\omega)$. If it is non-degenerate, then $(\mc{P}',\omega)$ is indeed a Cartan connection. This is equivalent with demanding that the tensor $\rP$ coming from $\nabla$ is non-degenerate.

We now turn to the inclusion $\mc{B} \subset \mc{A}$. The bundle $\mc{B}_0$ preserves both $H^*$ and (because it preserves split-Octonionic multiplication) $K$. Thus it must preserve the grading of $\T$, meaning that $\mc{B}_0 \subset \mc{A}_0$ -- or, equivalently, because of the algebra\"ic structures of these bundles, $\B_0 = \{\A_0,\A_0\}$. Decomposing $\A$ in terms of representations of $\B_0$, we get:
\be
\mc{A} = H' \oplus H \oplus (\mc{B}_0 \oplus \mc{L}) \oplus H^* \oplus (H')^*,
\ee
where $\mc{L}$ is the span of the grading section and we have defined $H' = T_{-2}$. Since $\mc{B}$ is of rank $14$, we can see that the projection $\B \to \mc{L}$ must be trivial. Since $\onab$ preserves $\B$, it must map $\B_0$ to $\B$. In particular, this means that $\nabla$ preserves $\mc{B}_0$ (thus preserves a volume form) and that for all sections $s$ of $\B_0$ and $X$ of $T$,
\beqa
\{X,s\} + \{\rP(X),s\} = -s(X) + s(\rP(X)).
\eeqa
is a section of $\B$. Since $\mc{B} \cap T^* = 0$, $s(\rP(X))$ must be determined by $s(X)$ -- consequently, the map $X \to \rP(X)$ is $\B_0$-invariant. In terms of representations of $\B_0$, $H \cong H^* \wedge H^* \cong (H')^*$, and $H^* \cong H'$ (these isomorphisms are given by the choice of any volume form on $H$). This means that $\rP|_H$ must be a multiple of the isomorphism $T_{-1} = H \cong (H')^* = T_2^*$, and $\rP|_{H'}$ must be a multiple of the isomorphism $T_{-2} = H' \cong H^* = T_1^*$. Thus $\rP_{11} = \rP_{22} = 0$ and we have three scenarios:
\begin{enumerate}
\item $\rP_{12} = \rP_{21} = 0$,
\item one of $\{\rP_{12},\rP_{21} \}$ in zero, the other is non-zero,
\item both $\rP_{12}$ and $\rP_{21}$ are non-zero.
\end{enumerate}
However, in the first case, $\B_0 \oplus T \oplus \rP(T)$ is not simple, and in the second, it isn't even an algebra bundle. In the third case, looking at the representations that $\rP$ identifies, we must have $\rP_{12} \propto \rP_{21}^t$. Then, considerations of the fact that $\B$ must be algebra\"ically closed and project trivially to $\mc{L}$ imply that $\rP_{12} = \rP_{21}^t$. Finally, since $\onab$ preserves $\B$, we must have 
\be
\nabla \rP = 0.
\ee
This is similar to, though not identical with, an Einstein involution (\cite{meein}). Since $\rP_{12}$ is an isomorphism, $\rP$ in particular is non-degenerate, so $(\mc{P}',\omega)$ is a Cartan connection.

Let $\sigma$ be the determinant of $\rP_{12}$; it is a volume form, and we will use this and $\mc{B}_0$ to give canonical identifications
\be
\A_{-2} \cong \A_1 \ \ \A_{-1} \cong A_2.
\ee
The reason for moving away from the tangent and cotangent space notations, is that we will be changing soldering forms when looking at $(\mc{P}',\omega)$. In this setting, the natural bundles are
\be
\A_{(k)}' = \sum_{j=-2}^k \A_j,
\ee
(in particular, $\A = \A'$). Using $\mf{p}'$ to define the filtration of $\mf{g}$ rather than $\mf{p}$, we can see that:
\be
\A_{(k)}' = \mc{G}_0 \times_{G_0} \mf{g}_{(k)}' = \mc{P}' \times_{P'} \mf{g}_{(k)}'.
\ee
This means that this splitting of $\A = \A'$ is a Weyl structure for $(\mc{P}',\omega)$, corresponding to the grading section $-E$. In any given splitting, the soldering form for a Tractor connection is given by sending any section $X$ of $T$ to:
\be
X \to \sum^{j=-1} \frac{1}{j}(\onab_X (-E))_j;
\ee
So for $(\mc{P}',\omega)$ the soldering form is:
\be
X \to \sum_{j=1} \frac{1}{-j}(\onab_X E)_j.
\ee
This identifies sections $X$ of $H$ with $\rP(X)$ of $\A_{-2}' = \A_{2}$ and sections $Z$ of $H'$ with $\rP(Z)$ of $\A_{-1}' = \A_1$. Thus $(\mc{P}',\omega)$ is a Cartan connection for the (free) $3$-distribution $H'$. We just need to show that it is normal.

Now $\onab$ is torsion free and has $G_2'$ holonomy, which means that the curvature $\kappa$ of $\onab$ is a section of $\wedge^2 T^* \otimes (\A_{(0)} \cap \B) = \wedge^2 T^* \otimes (\B_0)$. The harmonic curvature component of $\kappa$ (see Section \ref{har:cur}) is in $H \otimes H' \otimes \mc{A}'_0$, of homogeneity three. Other components of $\kappa$ must have higher homogeneity; the only possible candidate is $\kappa_{220}$, a section of $H' \wedge H' \otimes \B_0$. Since $\rP_{22} = 0$, this is precisely the $R^{\nabla}_{22}$, where $R^{\nabla}$ is the curvature of $\nabla$.

\begin{lemm}
$\kappa_{220} = R^{\nabla}_{22} = 0$.
\end{lemm}
\begin{lproof}
The Bianci identity for $\onab$ is $d^{\onab} \kappa = 0$, where $d^{\onab}$ is $\onab$ on $\mc{A}$ twisted with the exterior derivative $d$ on $\wedge^2 T^*$. Since both $\{,\}$ and $\rP$ are preserved by $\nabla$, $\nabla \kappa = \nabla R^{\nabla}$. Then for $X'$ and $Y'$ sections of $H'$ and with $Z$ a section of $H$,
\be
0 = (d^{\onab} \kappa)_{X',Y',Z}&=& (d^{\nabla} R^{\nabla})_{(X',Y',Z)} + \{X',\kappa_{(Y',Z)} \} + \{\rP(X'), \kappa_{(Y',Z)} \} + \textrm{ cyclic terms}\\ 
&=& 0 + \{X', (\kappa_{120})_{(Y',Z)} \} - \{Y', (\kappa_{120})_{(X',Z)} \} + \{\rP(X'), (\kappa_{120})_{(Y',Z)} \} \\
&& - \{\rP(Y'), (\kappa_{120})_{(X',Z)} \} + \{Z, (\kappa_{220})_{(X',Y')} \} + \{\rP(Z), (\kappa_{220})_{(X',Y')} \}
\ee
Now $\{\rP(Z'), (\kappa_{220})_{(X',Y')}$ is the only component taking values in $\mc{A}_{2}$, so it must vanish. Thus $\kappa_{220} = 0$.
\end{lproof}

Now we have $\kappa$ as a section of $H \otimes H' \otimes \B_{0}$. In particular $\kappa(H \wedge H) = 0$. Recall the definition of normality; that $\partial^* \kappa = 0$, where
\be
(\partial^* \kappa)(X) = \sum_{l} \{ Z^l, \kappa_{(Z_l, X)} \} - \frac{1}{2} \kappa_{(\{Z^l,X \}_{-}, Z_l)},
\ee
for $(Z_l)$ a frame for $T$ and $(Z^l)$ a dual frame for $T^*$. If we pick $\{Z_l\}$ so that it is the union of a frame of $H$ and a frame of $H'$, we know that $\{Z^l,X \}_{-} \wedge Z_l$ is zero or a section of $H \wedge H$ for all $X$. Thus the normality of $\kappa$ is entirely encapsulated in the fact that $\kappa$ is trace free:
\be
0= (\partial^* \kappa)(X) = \sum_{l} \{ Z^l, \kappa_{(Z_l, X)} \}.
\ee
We now need to calculate $(\partial^*)' \kappa$, where $(\partial^*)'$ is the operator for $(\mc{P}', \omega)$. Substituting in the soldering form for this second parabolic,
\be
((\partial^*)' \kappa)(X) = \sum_{l} \{ \rP(Z^l), \kappa_{(Z_l, X)} \} - \frac{1}{2} \kappa_{( \rP^{-1}(\{\rP^{-1}(Z^l), \rP(X) \}_{+}), Z_l)}.
\ee
Again, the second component vanishes, as $\kappa$ is zero on $H'\wedge H'$, leaving us with
\be
((\partial^*)' \kappa)(X) = \sum_{l} \{ \rP(Z^l), \kappa_{(Z_l, X)} \}
\ee
However, because $\rP$ is an isomorphism of $\B_0$-modules, and because $\kappa_{(Z_l, X)}$ is a section of $\B_0$, we have the equality: $\{ \rP(Z^l), \kappa_{(Z_l, X)} \} = \rP(\{ Z^l, \kappa_{(Z_l, X)} \})$. Thus
\be
((\partial^*)' \kappa)(X) = \sum_l \rP(\{ Z^l, \kappa_{(Z_l, X)} \}) = \rP(\partial^* \kappa)(X) = 0,
\ee
and $(\mc{P}',\omega)$ is normal.

Since the action of $\nabla$ on $\A_{1}$ and $\A_{2}$ is isomorphic with its action on $\A_{-2}$ and $\A_{-1}$, respectively, and since we may rewrite $\onab_X$ as
\be
\onab_X = \nabla_X + \rP^{-1}(\rP(X)), + (\rP(X)),
\ee
we can see that $\nabla$ is also the preferred connection for $(\mc{P}',\omega)$ in this splitting, and that its rho tensor is $\rP^{-1}$ conjugated by $\rP$ -- hence it is also $\rP$. Thus had we started with $(\mc{P}', \omega)$ and done the above construction, we would have ended up with $(\mc{P},\omega)$, and hence $(H')' = H$.
\end{proof}

\subsection{CR structures} \label{CR:Feff}
We aim to show here that there is a Fefferman construction for $\wh{G} = SO(4,3)$, $\wh{P}$ stabilises an isotropic $3$-plane, and $G = SO(4,2)$ while $P = (SO(2) \oplus GL(2)) \rtimes (\mbb{R}^2 \otimes \mbb{R}^{(2)}) \rtimes \wedge^2 \mbb{R}^2$ stabilises an isotropic $2$-plane.

Let $V = \mbb{R}^{(4,3)}$ and $B$ be an isotropic $3$-plane whose inclusion defines $\wh{P} \subset \wh{G}$. Let $W \cong \mbb{R}^{(4,2)}$ and fix an inclusion $W \subset V$ that defines $G \subset \wh{G}$.

Because of their signatures, $W$ and $B$ must be transverse, so their intersection $C = W \cap B$ is an isotropic $2$-plane. Defining $P$ as the stabiliser of $C$, it is evident that $G \cap \wh{P} \subset P$.

Now let $B'$ be the orthogonal projection of $B$ onto $W$. By construction, $C \subset B' \subset C^{\perp}$. The space $B'$ is equivalently defined by a line through the origin in $C^{\perp} / C$. The group $P$ acts via $SO(2)$ on this space of lines. Thus $G \cap \wh{P}$ lies as a codimension one subgroup in $P$. Then $\wh{G}$ is of dimension $21$, $\wh{P}$ of dimension $15$, $G$ also of dimension $15$, $P$ of dimension $10$ and $G \cap \wh{P}$ of dimension $9$. This implies that $G$ and $\wh{P}$ are transverse in $\wh{G}$, hence that the inclusion $\mf{g}/(\mf{g} \cap \wh{\mf{p}}) \to \wh{\mf{g}} / \wh{\mf{p}}$ is open. Thus we may do the Fefferman construction.

\begin{defi}[CR] A CR manifold is given by a contact distribution $K \subset TN$ with a complex structure $J$ on $K$. If $Q = TN/K$, and $q: TN \to Q$ is the obvious projection, there is a skew symmetric map $\mc{L}: K \times K \to Q$ given by $\mc{L}(X,Y) = q([X,Y])$ where $X$ and $Y$ are sections of $K$.

Integrability comes from using $J$ to split $K \otimes \mbb{C}$ as $K^{1,0} \oplus K^{0,1}$; the CR structure is integrable if $K^{0,1}$ is closed under the Lie bracket. This implies that $\mc{L}$ is of type $(1,1)$, that is that $\mc{L}(JX,JY) = \mc{L}(X,Y)$.
\end{defi}

\begin{theo}
The geometries modelled on $G/P$ are the 5 dimensional split signature CR geometries. If the CR structure is integrable and the Cartan connection is normal, the Cartan connection on the free $3$-distribution coming from the Fefferman construction is also normal.

Conversely, if the holonomy group of a normal Cartan connection for a free $3$-distribution reduces to $SO(4,2)$, it is locally the Fefferman construction over an integrable split signature CR geometry with normal Cartan connection.
\end{theo}
The rest of this section will be devoted to proving that theorem.

The first statement -- that the $G/P$ geometries are the CR geometries -- from the fact that the representation of $P$ as a parabolic is $\CR$, the same as for CR structures, combined with the following lemma:
\begin{lemm}
$Spin(4,2)_0 = SU(2,2)$.
\end{lemm}
\begin{lproof}
Consider the action of $SU(2,2)$ on $V = \mbb{C}^{(2,2)} \wedge \overline{\mbb{C}}^{(2,2)}$. $V$ carries a real structure on it from the action of the K{\"a}hler form, and a natural $(4,2)$ signature metric. Since $SU(2,2)$ is simple, and its action is non-trivial on this space, there is an inclusion
\be
\mf{su}(2,2) \hookrightarrow \mf{so}(4,2).
\ee
And then dimensional considerations imply that this is an equality. The maximal compact subgroup of $SO(4,2)_0$ is $S(O(4) \times O(2))_0$; the maximal compact subgroup of $SU(2,2)$ is $S(U(2) \times U(2))$. Then the Lemma is proved by looking at the fundamental groups of these Lie groups:
\be
\pi_1(SO(4,2)_0) &=& \mbb{Z}_2 \oplus \mbb{Z} \\
\pi_1(SU(2,2)) &=& \mbb{Z}. \\
\ee
\end{lproof}
Then it is easy to see that $P$ is the stabiliser of a complex null-line in $\mbb{C}^{(2,2)}$, demonstrating that these are split signature CR geometries (see \cite{CR}).

In order to demonstrate the normality conditions, we shall use both this Fefferman construction and the standard CR to conformal Fefferman construction (see \cite{CR}). Let $\wh{\wh{G}} = SO(4,4)$, with $\wh{\wh{P}}$ the stabiliser of a null line.  In details, if $(\mc{P}, \omega)$ is a split signature CR geometry, we have three related structures:
\be
(\ \mc{P}\ , \ \ \omega\ ) \ \ \ \ \ (\ \wh{\mc{P}}\ , \ \ \wh{\omega}\ ) \  \ \ \  \ (\ \wh{\wh{\mc{P}}}\ , \ \ \wh{\wh{\omega}}\ ),
\ee
where $\wh{\omega}$ is the Tractor connection for a free $3$-distribution while $\wh{\wh{\omega}}$ is a conformal Tractor connection.

We know that ${\wh{\omega}}$ is normal if and only if $\wh{\wh{\omega}}$ is normal (see \cite{bryskew}). But now consider the total inclusion of $G$ into $\wh{\wh{G}}$, given by composing the two Fefferman constructions.

If we complexify everything, we have $Spin(6) \subset Spin(7) \subset SO(8, \mbb{C})$. The spin representations of $Spin(6)$ are isomorphic with $\mbb{C}^4$, so decompose $\mbb{C}^8$ into two distinct components.

This implies that the action of $Spin(4,2) \subset Spin(4,3) \subset SO(4,4)$ on $\mbb{R}^{(4,4)}$ either decomposes it into two four dimensional components, or is irreducible on it (and preserves a complex structure on it). However, $SU(2,2) = Spin (4,2)$ does not have any four dimensional real representations (apart from the trivial one). Consequently the inclusion $SU(2,2) \subset SO(4,4)$ is the standard inclusion.

This means that the inclusion $G \subset \wh{\wh{G}}$ is the standard one. This generates $\wh{\wh{\omega}}$ via the iterated Fefferman construction. But this has to be the standard Fefferman construction. This implies that $\wh{\wh{\omega}}$ is normal if and only if $\omega$ is normal and the CR structure is integrable (see \cite{CR}, \cite{leitnersu} and \cite{leitnercom}). Consequently, $\wh{\omega}$ is normal if and only if $\omega$ is normal and the CR structure is integrable.

\begin{rem}
The inclusion $SU(2,2) \subset Spin(4,3)$ can be seen directly. $Spin (4,3)$ is defined as preserving a generic four-form $\lambda$ on $\mbb{R}^{(4,4)}$ (see \cite{baumspin}). $SU(2,2)$ on the other hand, preserves a K{\"a}hler form $\mu$, which can be seen as a section of $\wedge^2 V$ that is conjugate linear with respect to the volume form. It also preserve a complex volume form $v \in \Gamma(\wedge^{(4,0)} V_{\mbb{C}})$. The inclusion of $SU(2,2)$ into $Spin(3,4)$ is given by the generic four form:
\beqa \label{form:su}
Re(v) - (\mu)^2.
\eeqa
\end{rem}

Now we need to show the converse. Let $(M, \wh{\mc{P}}, \onab)$ be a normal Cartan connection for a free $3$-distribution. Assume the holonomy group of $\onab$ reduces to $SO(4,2)$ -- equivalently, that there is a section $\tau$ of $\mc{T}$, of negative norm squared, such that $\onab \tau = 0$. Define $R \in \Gamma(H)$ as $\pi^2(\tau)$; since $\tau$ is of negative norm squared, $R$ is never-zero.

Define the bundle $\wt{K}$ as $H \oplus [H,R]$. It is a bundle of rank five. Let $N$ be the manifold got from $M$ by projecting along the flow of $R$.
\begin{prop}
$N$ carries a CR structure, and the contact distribution $K$ in $TN$ is the projection of $\wt{K}$ to $N$. The complex structure $J$ on $K$ is given by the action of $R$.
\end{prop}

\begin{proof}
We first need to show that $[\wt{K}, R] = \wt{K}$. Let $L$ be the line subbundle of $\mc{T}$ generated by $\tau$. Pick any preferred connection $\nabla$ such that in the splitting it defines, $\tau = (\alpha,0,R)$ and $\alpha(R) = 1$ (to show this is possible, use any preferred connection to get a splitting $\tau = (\alpha',f,R)$ and change the splitting by the action of an $\Upsilon$ where $\Upsilon \llcorner R = f$). Call such a $\nabla$ an $L$-preferred connection. Now, since $\onab \tau = 0$, we can see that $\nabla_X \alpha = \nabla_X R = 0$ for any section $X$ of $H$ while $\nabla_U R = - \{U,\alpha\}$ for $U$ a section of $T_{-2}$ (we negelect the action of $\rP(X)$ and $\rP(U)$ on $(0,0,R)$ as this action is trivial).

We may choose sections $X$ and $Y$ of $H$ that obey the following properties:
\begin{enumerate}
\item $X$, $Y$ and $R$ form a frame of $H$,
\item $\nabla_R X = \nabla_R Y = 0$,
\item $\alpha(X) = \alpha(Y) = 0$, 
\end{enumerate}
(for instance, we could define $X$ and $Y$ obeying the algebraic properties on a submanifold transverse to $R$, and extend by parallel transport along $R$; then the relation $\nabla_R R = \nabla_R \alpha = 0$ ensures the algebraic properties are preserved). Since $\onab$ is torsion-free,
\be
[R, X] &=& \nabla_R X - \nabla_X R - \{R,X\} \\
&=& \{R,X\}.
\ee
Similarly,
\be
[R,[R,X]] &=& \nabla_R \{R,X\} - \nabla_{\{R,X\}} R - 0 \\
&=&  \{\{R,X\}, \alpha \} = -X.
\ee
The same hold for $Y$ and $\{R,Y\}$. Consequently $[R, \wt{K}] = \wt{K}$ and $\wt{K}$ projects to a distribution $K$ in $N = M/R$. This distribution must be a contact distribution, by the properties of the Lie bracket on $\wt{K}$. Let $r$ be any coordinate on $M$ such that $R\cdot r = 1$. Then the vector fields
\beqa \label{sec:CR}
\cos(r)X - \sin(r)\{R,X\},& \sin(r)X - \cos(r)\{R,X\} \\
\cos(r)Y - \sin(r)\{R,Y\},& \sin(r)Y - \cos(r)\{R,Y\}
\eeqa
are $R$-invariant, hence lifts of vector fields in $K$. Since we have this explicit form, we can see that the Lie bracket of $X$ and $Y$ with $R$ generates an endomorphism of $\wt{K}$ that descends to an automorphism $J$ of $K$, squaring to minus the identity.

Changing to another $L$-preferred connection will change $X$ and $Y$ by adding multiples of $R$. This will change neither their projections nor the properties of $J$. This means that the CR structure is well defined.

And $M$ must be the Fefferman construction over this CR structure. This implies that CR structure must be integrable and that $\onab$ descends to a normal CR Tractor connection on $N$.
\end{proof}

One interesting consideration: there are non-trivial morphisms of the $3$-distribution that cover the identity on the underlying CR structure. Any diffeomorphism $\phi: M \to M$ generated by a flow on $R$ will change the distribution $H$, but since $\phi$ projects to the identity on $N$, it leaves the underlying CR structure invariant.

The forgoing means that all the results on CR holonomy (equivalently, conformal holonomy contained in $\mf{su}(2,2)$) have equivalent formulations in terms of free $3$-distributions. See papers \cite{leitnersu}, \cite{CR} and \cite{leitnercom}; paper \cite{mecon} has some Einstein examples as well. This implies, for instance, that holonomy reductions to $SU(2,1)$ exist (whenever $N$ is a Sasaki-Einstein manifold with the correct signature and sign of the Einstein coefficient). From the free $3$-distribution point of view, this corresponds to a complex structure on the complement of $\tau$ in $\mc{T}$.

Similar consideration exist for a holonomy reduction to $Sp(1,1) \subset SU(2,2)$, with the quaternionic analogue of CR spaces.

\subsection{Lagrangian contact structures}
Lagrangian contact structures (see for example \cite{complexconf}) geometries generated by another real form of the parabolic that models CR structures.

\begin{defi}
A Lagrangian contact structure is given by a contact distribution $K$ on a manifold of dimension $2m+1$, together with two bundles $E$ and $F$ of rank $m$ such that $K = E \oplus F$, $[E,E] \subset K$ and $[F,F] \subset K$.
\end{defi}
The structure is integrable if both $E$ and $F$ are integrable. The parabolics are given by $G = SL(m, \mbb{R})$ while $P = (\mbb{R} \oplus GL(m-2, \mbb{R}) \rtimes (\mbb{R}^{m} \oplus \mbb{R}^{m*}) \rtimes \mbb{R}$. Then we have:
\begin{lemm}
$Spin(3,3)_0 = SL(4,\mbb{R})$
\end{lemm}
\begin{lproof}
Consider the action of $SL(4,\mbb{R})$ on $V = \wedge^2 R^4$. Since $SL(4,\mbb{R})$ preserves a volume form which is an element of $\wedge^4 \mbb{R}^{4*} \subset \odot^2 (\wedge^2 \mbb{R}^4)^*$, it preserves a metric on $V$, of split signature. Then since $SL(4,\mbb{R})$ is simple and acts non-trivially, we get an algebra inclusion $\mf{sl}(4,\mbb{R}) \subset \mf{so}(3,3)$ and the dimensions imply equality.

The maximum compact subgroup of $SL(4,\mbb{R})$ is $SO(4)$ while the maximum compact subgroup of $SO(3,3)_0$ is $S(O(3) \times O(3))_0$. Consequently the result flows from:
\be
\pi_1 (SL(4,\mbb{R})) = \mbb{Z}_2, \\
\pi_1 (SO(3,3)_0) = \mbb{Z}_2 \times \mbb{Z}_2.
\ee
\end{lproof}

Given this, the results for CR structures go through almost verbatim to this new setting, bar one subtlety: $\mbb{R}^{(3,3)}$ need not be transverse to a given isotropic $3$-plane in $\mbb{R}^{(4,3)}$. So we may need to restrict our results to open dense subsets of our manifolds. The inclusion $GL(4,\mbb{R}) \subset Spin(4,3) \subset SO(4,4)$ is the standard inclusion that decomposes $\mbb{R}^{(4,4)}$ as $\mbb{R}^4 \oplus \mbb{R}^{4*}$. Summarising all these results:
\begin{theo}
Let $N$ be a five dimensional integrable Lagrangian contact manifold. Then there is a Fefferman construction for $N$ to a free $3$-distribution on a manifold $M$, the total space of a circle or line bundle over $N$. The Tractor connection on $M$ is normal if and only if the Tractor connection on $N$ is normal.

Conversely, if $M$ is a free $3$-distribution geometry with normal Tractor connection $\onab$, and the holonomy group of $\onab$ reduces to $SO(3,3)$, then \emph{on an open dense set}, $M$ is locally the Fefferman space of a integrable, normal Lagrangian contact manifold.
\end{theo}

\vspace{5mm}
\begin{center}
	{\Large{{\bf APPENDIX}}}
\end{center}
\vspace{-5mm}
\appendix
\section{Higher dimensional holonomy}
We can look briefly at the higher dimensional free distribution, where $H$ is of rank $n > 3$ and the whole manifold is of dimension $n(n+1)/2$. Considerations of the harmonic curvature (see Section \ref{har:cur}) imply that the Tractor connection is torsion-free if and only if it is flat. Also, the Fefferman construction of Section \ref{spin} gives an almost-spinorial structure, not a conformal one. So we can expect that the situations in higher dimensions is not a simple extention of that in rank $3$.

Indeed, there are some fascinating potential structures on these manifolds -- sub-Riemannian structures, restrictions to lower dimensional distributions, and special classes of preferred connections for each preserved subbundle of the tractor bundle. These are detailled in an old arXiv paper of mine, \cite{meskew}. Existence results for the normal, non-flat cases are generally lacking, however. With one major exception:
\begin{theo}
For every $n$, and every $p \leq \frac{n(n-1)}{2} - 3$, there exists a manifold $M$ of dimension $n(n+1)/2$, with $H$ a free distribution of rank $n$ on it, such that the normal Tractor connection $\onab$ generated by $H$ has holonomy algebras isomorphic to $\mbb{R}^p$.
\end{theo}
\begin{proof}
Notice that for $n=3$, we require $p\leq 0$, so the result is only relevant for $n>3$. Now the homogeneous model of a free distribution is given locally (\cite{meskew}) by a frame $\mc{F} = \{X_i,Y_{j|k} \}_{j<k}$ of $T$, with $i,j,k$ running from $1$ to $n$. All these vector fields commute, with the exception of:
\be
[X_j, X_k] = Y_{j|k},
\ee
for $j <k$. The distribution $H$ is identified with the span of the $X_i$'s. If we pick local coordinates $x_i$ and $y_{j|k}$ on $M$, we may define the frame via:
\be
Y_{j|k} = \frac{\partial}{\partial y_{j|k}} \ \ \textrm{and} \ \ X_i = \frac{\partial}{\partial x_{i}} - \sum_{p=i+1}^n x_p Y_{i|p} 
\ee

We shall modify this construction slightly, to get the desired result. We will change the frame $\mc{F}$ by replacing with $X_1$ with $X_1' = X_1 + y_{1|2} Y_{3|4}$. The Lie brackets for this new frame are all the same as for the old frame, with the exception of
\be
[Y_{1|2}, X_1'] = Y_{3|4}.
\ee
Now define $\nabla$ to be the flat connection annihilating all the vector fields in $\mc{F}$, $H$ as the span of $X_1'$ and the $X_i$'s ($i \neq 1$), and $T_{-2}$ as the span of the $Y_{j|k}$'s. Let $Z$ be any vector field on $M$. We wish to show that the Tractor connection defined as $\onab_Z = Z + \nabla_Z$ is normal. Its curvature $\kappa$ is trivially flat, apart from the term
\be
\kappa(Y_{1|2}, X_1') = Tor^{\nabla}(Y_{1|2}, X_1') = Y_{3|4}.
\ee
Looking at the formula (\ref{normal:form}) for the normality condition on $\partial^* \kappa$, we can see that the only terms in which $\kappa(Y_{1|2}, X_1')$ appear are:
\be
\partial^* \kappa(X_1') &=& \{(Y_{1|2})^*, \kappa(Y_{1|2}, X_1') \} = \{(Y_{1|2})^*,Y_{3|4}\} = 0, \\
\partial^* \kappa(Y_{3|4}) &=& \{(X_1')^*, \kappa(X_1', Y_{1|2}) \} = - \{(X_1')^*,Y_{3|4}\} = 0,
\ee
so $\onab$ is normal. Note that in this case, the full curvature is isomorphic to the harmonic curvature. Now $\nabla \kappa = 0$ as $\nabla$ is flat, and the remaining term in $\onab \kappa$ vanishes for reasons of homogeneity. This implies that $\onab_{A_1} \onab_{A_2} \ldots \onab_{A_{l}} \kappa(B,C)$ is a $C^{\infty}(M)$-multiple of $Y_{3,4}$ for all sections $A_1, \ldots A_l, B$ and $C$ of $TM$. So considerations of infinitesimal holonomy (see \cite{KobNa}) tell us that the holonomy algebra of $\onab$ is the span of $Y_{3|4}$ (in the splitting of $\mc{A}$ defined by $\nabla$).

Similarly, we may substitute for $X_1$ and $X_2$ via
\be
X_1' &=& X_1 + \sum_{j=3,k>j}^n \frac{1}{(\beta(j,k))!}(y_{1|2})^{\beta(j,k)} Y_{j|k} + \sum_{j=4} \frac{1}{(\gamma(j))!} (y_{1|3})^{\gamma(j)} Y_{2|j} \\
X_2' &=& X_2 + \sum_{j=4}^n \frac{1}{(\delta(j))!} (y_{2|3})^{\delta(j)} Y_{1|j},
\ee
where $\beta$, $\gamma$ and $\delta$ are any \emph{injective} functions to the strictly positive integers.

Defining $\nabla$ flat as before and $\onab_Z = Z + \nabla_Z$, the above argument demonstrates that $\kappa$ is normal, since the only terms in $\kappa$ are multiplies of $Y_{j|k}^* \otimes X_l^* \otimes Y_{a|b}$ with $a$ and $b$ distinct from $j$, $k$ and $l$. Then if we define $(\onab_Z)^l$ as $\onab_Z$ iterated $l$ times, one has, at the origin:
\be
(\onab_{Y_{1|2}})^{\beta(j,k)-1}  \big(\kappa(X_1',Y_{1|2})\big) & = Y_{j|k} +O(1) & \ \ \textrm{for $k > j > 2$}, \\
(\onab_{Y_{1|3}})^{\gamma(j) - 1} \big(\kappa(X_1',Y_{1|3})\big) & = Y_{2|j} + O(1) & \ \ \textrm{for $ j > 3$},\\
(\onab_{Y_{2|3}})^{\delta(j) - 1} \big(\kappa(X_1',Y_{2|3})\big) & = Y_{1|j} + O(1) & \ \ \textrm{for $j > 3$}.
\ee
The only frame elements that cannot be generated in this way are $Y_{1|2}$, $Y_{1|3}$ and $Y_{2|3}$; in fact, they cannot appear anywhere in the infinitesimal holonomy. Thus we have a holonomy algebra isomorphic with $\mbb{R}^q$ where $q = n(n-1)/2 - 3$. Smaller holonomy algebras can then be generated by removing terms from the above expressions for $X_1'$ and $X_2'$.
\end{proof}

\bibliographystyle{amsalpha}
\bibliography{ref}

\end{document}